\documentclass[12pt]{article}

\usepackage{amsmath,amsfonts,amssymb}

\title{Wave packets in the Schwartz space of  a reductive $p$-adic symmetric space}

\author{ Patrick Delorme\thanks{ P. Delorme is a member of the Institut Universitaire de France} , Pascale Harinck}

\date{}

\def\a{\mathfrak{a}}

\def\b{\mathfrak{b}}

\def\si{\sigma}
\def\ep{\varepsilon}

\def\A{\mathcal{A}}
\def\CC{\mathcal{C}}

\def\O{\cal O}

\def\LL{\Lambda}
\def\l{\lambda}
\def\aa{\alpha}

\def\DD{\Delta}

\def\F{\mathbf{F}}

\def\mm{\mathcal}

\def\N{\mathbb{N}}

\def\Z{\mathbb{Z}}
\def\R{\mathbb{R}}
\def\C{\mathbb{C}}

\def\fd {\hspace{0.35cm} \raise
-0.5mm\hbox{$\blacksquare$}\\}

\def\qed{{\null\hfill\ \raise3pt\hbox{\framebox[0.1in]{}}\break\null}}
\newtheorem{Theorem}{Theorem}[section]
\newtheorem{Proposition}[Theorem]{Proposition}
\newtheorem{Lemma}[Theorem]{Lemma}

\newtheorem{Remark}[Theorem]{Remark}
\newtheorem{Definition}[Theorem]{Definition}

\def\ste{\par\smallskip\noindent}

\def\dem{ {\em Proof~: \ste }}

\def\beq{\begin{equation}}
\def\eeq{\end{equation}}

\def\bb{\backslash}

\def\A {\cal A}

\def\O {\mathcal O}

\def\DD{\mathcal{D}}
\newenvironment{res}
               {\begin{equation}\begin{minipage}{0.85\textwidth}}
               {\end{minipage}\end{equation}}
\def\ber{\begin{res}}
\def\eer{\end{res}}

\def\un{ \underline}
\def\mm{\mathcal}

\catcode`\Ž=\active\def Ž{\'e}
\catcode`\=\active\def {\`e}
\catcode`\ˆ=\active\def ˆ{\`a}
\catcode`\=\active\def {\`u}
\catcode`\=\active\def {\^{e}}
\catcode`\"=\active\def "{\^{i}}
\catcode`\™=\active\def ™{\^{o}}
\catcode`\ž=\active\def ž{\^{u}}
\catcode`\=\active\def {\c{c}}

\begin{document}

% Insert title

\maketitle

\begin{abstract}
We form wave packets in the Schwartz space of a reductive $p$-adic symmetric space for certain famillies of tempered functions. We show how to construct such families from Eisenstein integrals.
\end{abstract}

\medskip


\section{Introduction}
Let $G$ be the group of $\F$-points of an algebraic group, $\un{G}$, defined over $\F$, where $\F$ is a nonarchimedean local field of characteristic different from 2.
Let $H$ be the group of $\F$-points of an open $\F$-subgroup of the fixed point group of a rational involution of $\un{G}$ defined over $\F$.

We introduce the space $\mm{A}_{temp}(H\bb G)$ of smooth tempered  functions on $H\bb G$. They are the tempered functions  which are  generalized coefficients of an $H$-fixed linear form $\xi$ on an admissible  $G$-module $V$, when $V$ and $\xi$ varies.

 Using the theory of the constant term (cf. \cite{L}, \cite{KT1}), we introduce the weak constant term of elements of $\mm{A}_{temp}(H\bb G)$ as it was made in \cite{W} for tempered functions on  the group.
 
Then, we introduce families of elements of $\mm{A}_{temp}(H\bb G)$ of type I, by conditions on their exponents. Then the conditions are strengthened  to introduce families of type I', and we add  conditions on the weak constant term to define families of type  II'. This is the analogue of  the families used in \cite{BaCD} for the real case. 

Important examples  of such families are given (cf. Theorem \ref{eisII'}) in terms of Eisenstein integrals, due to the main results of \cite{CD}. 

Then, following \cite{BaCD} for the real case, which was largely inspired by the work of Harish-Chandra \cite{H-C}, and   \cite{W}, we show that one can form wave packets in the Schwartz space for such families (Theorem \ref{theo1}). Notice also that  the intermediate Proposition \ref{fq+} is the analogue of the important  Lemma 7.1 of \cite{A}. 
 
The recent work of Sakerallidis and Venkatesh \cite{SaV} on spherical varieties includes in particular the $L^2$-Plancherel formula for $H\bb G$, when $G$ is split and the characteristic of $\F$ is equal to zero. It should be possible using our  result  to determine the Fourier transform of the Schwartz space for these symmetric spaces. This  should be entirely analogous to the work [DO] for affine Hecke algebras.\ste

\noindent {\bf Acknowledgments.} We thank warmly the referee for his very pertinent mathematical comments and  his careful remarks on our presentation. We thank also Omer Offen for his collaboration when this work was intended to be used for truncation on  some particular symmetric spaces.

\section{The map $H_G$ and the real  functions $\Theta_G, \Vert . \Vert $ and $N_d$ on $H\bb G$}

\noindent{\bf 2.1. Notation.}\ste
\setcounter{equation}{0}
If $E$ is a vector space, $E'$ will denote its dual. If $E$ is real, 
$E_\C$ will denote its complexification.\ste If $G$ is a group,
$g\in G$ and  $X$ is a subset of $G$, $g.X$ will denote $gXg^{-1}$. If $J$ is a subgroup of $G$, $g\in G$
and $(\pi, V)$ is a representation of $J$, $V^J$ will denote the space of
invariant elements of $V$ under $J$ and $(g\pi, gV)$ will denote the representation of $g.J$ on $gV:=V$ defined by: $$(g\pi)(gxg^{-1}): = \pi(x), x\in J.$$ We will denote by $(\pi', V')$ the dual representation of $G$ in the algebraic dual vector space $V'$ of $V$. \\If $V$ is a vector space of vector valued functions  on $G$ which are invariant by right (resp., left ) translations, we will denote by $\rho$ (resp., $\l$) the right (resp., left) regular representation of $G$ in $V$.\\ If $G$ is locally compact, $d_lg$ will denote a left invariant Haar measure on $G$ and $\delta_G$ will denote the modulus function. \\ Let  $\F$
be a non archimedean local field with finite residue field $\F_q$. Unless specified we assume: 
\ber \label{22}  The characteristic of $\F$ is different from 2. \eer
Let $\vert.\vert_{\F} $ be the  normalized absolute value  of $\F$.\\
We will use conventions like in \cite{W}.
 One considers various algebraic groups defined over $\F$,  and a sentence
like:  
 \ber \label{corse} " let $A$ be a split torus
 "  will mean '' let 
$A$ be the  group of $\F$-points   of a torus,  $\underline A$,   defined and split  over $\F$ ".\eer 
With these  conventions, let $G$ be a connected reductive linear algebraic
group. Let $\tilde{A}_G$ be the maximal split torus of the center of $G$. The change to standard notation will become clear later.
\ste 
Let $A$ be a split  torus of $G$. Let  $X_{*}(A)$ be the group  of one-parameter subgroups of  
 $A$. This is a free abelian group of finite type. Such a group will be
called a lattice. One fixes a uniformizer 
 $\varpi$ of  $\F$. One denotes by  $\Lambda (A)$ the image of  $X_{*}(A)$ in 
$A$ by  the morphism of groups  $\underline \lambda \mapsto \underline
\lambda (\varpi)$. By this morphism $\LL(A)$ is isomorphic to   $X_{*}(A)$. \ste 
If  $J$ is an algebraic group, one denotes by $\rm{Rat} (J)$ 
the  group of its rational characters defined over $\F$.   Let us define:
$$\tilde{\a}_{G}= \rm{Hom}_{\Z} (\rm{Rat} (G),\R).$$ The restriction of rational characters 
 from  $G$ to  $\tilde{A}_{G}$  induces an  isomorphism: 
\beq \label{iso} \rm{Rat}(G)\otimes _{\Z}\R \simeq \rm{Rat}(\tilde{A}_{G})\otimes _{\Z} \R.\eeq
Notice that $\rm{Rat}(\tilde{A}_G)$ appears as a  generating lattice in the dual space $\tilde{\a}'_G$ of $\tilde{\a}_G$ and: \beq \label{aotimes} \tilde{\a}_G' \simeq \rm{Rat}(G) \otimes _\Z \R. \eeq   
One has the canonical map  $\tilde{H}_{G}: G \rightarrow \tilde{\a}_{G}$
 which is defined by: 
\beq \label{H} e^{\langle \tilde{H}_{G}(x), \chi\rangle}= \vert \chi (x)\vert_{\F}, \> x\in G,
\chi
\in \rm{Rat} (G).\eeq The kernel of  
$\tilde{H}_{G}$, which is denoted by $\tilde{G}^1$,   is the intersection of the kernels   of $\vert \chi \vert_{\F}$ for all  character $\chi \in \rm{Rat}
(G)$  of
$G$. One defines
$X(G)= \rm{Hom} (G/\tilde{G}^{1}, {\C}^*)$, which is the group of unramified characters of $G$.  One will
use similar notation for Levi subgroups of $G$.\ste  
   One denotes by  $\tilde{\a}_{G, \F}$ (resp.,  $\tilde{\tilde{\a}
}_{G, \F}$)  the image of  $G$ (resp.,  $\tilde{A}_{G}$) by  $\tilde{H}_{G}$. Then 
$G/\tilde{G}^{1}$ is  isomorphic to the lattice 
$\tilde{\a}_{G,\F}$.\ste 
There is a surjective map: 
\beq  \label{surjection}(\tilde{\a}'_{G})_{\C}\rightarrow X(G)\rightarrow
1 \eeq   denoted by $ \nu \mapsto \chi_\nu$ which associates to 
$\chi\otimes s$, with $\chi  \in \rm{Rat}(G) $, $s\in \C$,    the character  $g\mapsto \vert \chi (g)\vert_\F 
^{s}$ (cf. \cite{W}, I.1.(1)).  In other words: 
\beq  \label{cnu}\chi_\nu (g) =   e^{ \langle \nu, \tilde{H}_G (g)\rangle}, g \in G, \nu \in(\tilde{\a}'_{G})_{\C}.\eeq
 The kernel is a lattice and it defines a structure of a
complex algebraic variety on  
$X(G)$ of dimension   
$dim_{\R}\tilde{\a}_{G}$.  Moreover $X(G)$ is an abelian  complex Lie group whose Lie algebra is equal to $(\tilde{\a}'_{G})_{\C}$. \ste
 If  $\chi$ is an element of $ X(G)$, let  $\nu$ be an element of $
\tilde{\a}_{G,\C}'$ such that $\chi_\nu= \chi $.  The real part  $\rm{Re}\>\nu \in \tilde{\a}'_{G}$ is independent from the choice of $\nu$.  We will denote it  by $\rm{Re} \> \chi$. If  $\chi \in \rm{Hom} (G,
\C^{*}) $ is continuous, the character $\vert \chi \vert$ of $G$ is an element of  $X(G)$.  One sets $\rm{Re}\> \chi=  \rm{Re}\> \vert \chi\vert $. Similarly,  if  $\chi \in \rm{Hom}(\tilde{A}_{G},
\C^{*})$ is continuous, the character $\vert \chi \vert$ of $\tilde{A}_G$ extends uniquely  to  an element of  $X(G)$  with values in  $\R^{*+}$, that we will denote again by 
  $\vert
\chi \vert$ and one sets  $\rm{Re} \> \chi=  \rm{Re}\> \vert \chi
\vert$.\\
 If $P$ is a parabolic subgroup of $G$ with Levi subgroup $M$,
we  keep the same notation with $M$ instead of $G$.\\ The inclusions
$\tilde{A}_{G}\subset \tilde{A}_{M}\subset M\subset G$  determine a surjective  morphism    $\tilde{\a}_{M, \F}\rightarrow  \tilde{\a}_{G, \F}$ (resp., an injective  morphism, 
 $\tilde{ \tilde 
{\a}}_{G, \F}   \rightarrow\tilde{ \tilde 
{\a}}_{M, \F}$) which extends uniquely to a surjective linear map between  $\tilde{\a}_{M}$ and  $\tilde{\a}_{G}$, (resp., injective map  between
$\tilde{\a}_{G}$ and 
$\tilde{\a}_{M}$).
The second map allows to identify  $\tilde{\a}_{G}$ with a subspace of 
$\tilde{\a}_{M}$ and the kernel of the first one, $\tilde{\a}^{G}_{M}$,  satisfies: 
\ber $$\label{oplus} \tilde{\a}_{M}= \tilde{\a}^{G}_{M}\oplus \tilde{\a}_{G}.$$
\eer
\\
If an unramified character of $G$ is trivial on $M$, it is trivial on  any maximal compact subgroup of $G$ and on the unipotent radical of $P$, hence on $G$. It allows to identify $X(G)$ with a subgroup of $X(M)$. Then $X(G)$ is the analytic subgroup of $X(M) $ with Lie algebra $(\tilde{\a}'_G)_\C \subset (\tilde{\a}'_M)_\C$. This follows easily from  (\ref{cnu}) and (\ref{oplus}). 

  One has (cf.  \cite{D2}, (4.5)), \ber \label{45} The map $\LL(\tilde{A}_G )\to G/\tilde{G}^1$ is injective and allows to identify $\Lambda(\tilde{A}_G)$ with the subgroup $\tilde{H}_{G}(\tilde{A}_G)$ of $\tilde{\a}_G$. \eer  
Let ${\underline G}$ be the algebraic group defined over $\F$ whose group of $\F$-points is $G$. Let $\si$ be a rational involution of ${\underline G}$  defined over
$\F$.  Let 
$H$ be  the group of $\F$-points of an open $\F$-subgroup  of the fixed
point set of $\si$.  We will also denote by  $\sigma$  the restriction of $\si$ to
$G$.\\A split torus $A$ of $G$ is said to be  $\si$-split if $A$ is contained in  the set of elements of $G$ which are antiinvariant by $\si$.  Now we explain the change to standard  notation: $A_G$  will denote the maximal $\si$-split torus of the center of $G$.\\
Let $\tilde{A}$ be a $\si$-stable  split torus of $G$. The involution   $\si$ induces an involution, denoted in the same way,  on  $\tilde{\a}:=\tilde{\a}_{\tilde{A}}$. Let $\tilde{A}^\si$ (resp., $\tilde{A}_\si$) be the  maximal split torus in  the group of elements of $\tilde{A}$ which are  invariant (resp., antiinvariant) by $\si$.  Then $\tilde{\a}^\si$ (resp., $\tilde{\a}_\si$) is identified   with the set of invariant   (resp., antiinvariant) of $
\tilde{\a}$ by $\si $    and $\tilde{A}_\si$ is the maximal $\si$-split torus of $\tilde{A}$. \\ In particular,  
 one has $A_G=(\tilde{A}_G)_\si$ and $\tilde{\a}_G= \tilde{\a}_G^\si\oplus \a_G$ where $\tilde{\a}_G^{\si}$ (resp., $\a_{G}$) is the space of invariant (resp., antiinvariant) elements of $\tilde{\a}_G $ by $\si$.

We define  a morphism of groups $H_{G} : G \to  \a_{G} $ which is the composition of $\tilde{H}_G$ with the projection on $ \a_{G}$ parallel to $\tilde{\a}_G^\si$.  We remark that, as is  seen easily,  $\tilde{H}_G$ commutes with $\si$. Hence $H_{G} $ is trivial on $H$. \\ We denote by $G^1$ the kernel of $H_G$, which contains $H$. It contains also $\tilde{G^1}$, hence it is open in $G$. We denote by $\a_{G,\F}$ the image of $H_G$.
Let  $X(G)_\si$ be the connected component  of the  group of antiinvariant elements of $X(G)$. Then  $X(G)_\si$ is the analytic subgroup of $X(G)$ with Lie algebra  $(\a'_G)_\C \subset (\tilde{\a}'_G)_\C $. The elements of $X(G)_\si$ are precisely the characters of $G$ of the form $$\chi_\nu(g)= e^{\langle \nu, H_G(g)\rangle }, \nu \in (\a'_G)_\C, g\in G.$$ They are exactly the characters of the lattice $G^1\bb G$. The group $X(G)_\si$ has a natural structure of complex algebraic group. We denote  by $X(G)_{\si, u}$ the group of unitary elements of $X(G)_\si$.\\
 One has
\ber \label{llam} The group  $\Lambda(A_G)$ is identified by $H_G$ with  $H_G(A_G)$.\eer
Let  $\tilde{A}$  be a  maximal  split torus of $G$. Let $M$ be the centralizer of $\tilde{A}$ in $G$. Let us show the following assertion.
\ber \label{mweight}$\tilde{H}_M (\tilde{A})$  contains a   multiple by $k \in \R^{+*}$ of the coweight lattice of the root system $\Sigma\subset (\tilde{\a}_M^G)'$  of $\tilde{A}$ in the Lie algebra of $G$.
 Here the coweight lattice is the dual  lattice in $\tilde{\a}_M^G$  of the root lattice.\eer 
It is clear that it suffices to prove the assertion for one maximal split torus. Let $\tilde{A}'$ be a maximal split torus of the derived group $G_{der}$ of $G$. Let $\underline{\tilde{A}}= \underline{\tilde{A}'} \underline{\tilde{A}_G}$. This is a maximal $\F$-split torus of $\underline{G}$ for reasons of dimension. The intersection $F$ of $\underline{\tilde{A}'} $ and $\underline{\tilde{A}_G}$ is finite. Hence one has the exact sequence  $$1\to F \to  \underline{\tilde{A}'}\times  \underline{\tilde{A}_G} \to \underline{\tilde{A}} \to 1. $$ 
 Going to $\F$-points,  the long exact sequence in cohomology implies that 
$  \tilde{A}\tilde{A}_G'$ is of finite index in $\tilde{A}$. Hence the image of $\tilde{A}'\tilde{A}_G$ by $\tilde{H}_{M}$ is of finite index in the image of $\tilde{A}$.
The image of $\tilde{A}'$ (resp., $\tilde{A}_G$) in $\tilde{\a}_M$  by $\tilde{H}_M$ is contained in $\tilde{\a}_M^G$ (resp., $\tilde{\a}_G$) and is a lattice $\LL_1$ (resp., $\LL_G$) generating $\tilde{\a}_M^G$ because $\LL_1+ \LL_G$ is of finite index in $\LL=\tilde{H}_{M} (\tilde{A})$ which generates $\tilde{\a}_M$. Hence the rank of $\LL_1$ is equal to the dimension of 
$\tilde{\a}_M^G$. 
 The values of the normalized  absolute value of $\F$ are of the form $q^n, n \in \Z$. From the definition of $\tilde{H}_M$,  one sees that $\LL_1$ is included in $(\log q)\LL_2$  where $\LL_2\subset \tilde{\a}_M^G $ is the coweight   lattice  of  $\Sigma$.
Both are lattices of the same rank, for reasons of dimension.
Our claim follows from the following assertion:
\ber Let $\LL_1\subset \LL_2$ be two lattices of the same rank. Then there exists $n \in \N^*$ such that  $n \LL_2\subset \LL_1$, \eer
which follows by inverting the matrix, with integral entries, expressing a basis of $\LL_1$ in a basis of $\LL_2$.\\
Let $A$ be a maximal $\si$-split torus of $G$ and let $\tilde{A}$ be a $\si$-stable maximal split torus of $G$ which contains $A$. The roots of $A$ in the Lie algebra of $G$ form a root system (cf. \cite{HW}, Proposition 5.9).  Let $M$ be the centralizer in $G$ of $A$, which is $\si$-invariant. One has $A=A_M$. One deduces like  (\ref{mweight}) that: 
 \ber \label{Mweight}$\LL(A) \subset \a$ contains  a   multiple by $k \in \R^{+*}$ of the coweight lattice  of the root system of $A$ in the Lie algebra of $G$. \eer
 A  parabolic subgroup $P$ of $G$  is called a $\sigma$-parabolic subgroup if  $P$ and
$\sigma(P)$ are opposite  parabolic subgroups. Then $M:=P\cap \sigma(P)$ is the  $\sigma$-stable
Levi subgroup of $P$.    If $P$ is such a parabolic subgroup, $P^-$ will denote $\si(P)$. \\   The sentence :
 ''Let $P=MU$ be a parabolic subgroup of $G$'' will mean that $U$ is the unipotent radical of $P$ and $M$ a Levi subgroup of $G$. If moreover $P$ is a $\si$-parabolic subgroup of $G$, $M$ will denote its $\si$-stable Levi subgroup.\\
Let $P=MU$ be a $\si$-parabolic subgroup of $G$.   Recall that $A_M$ is  the maximal $\si$-split torus of the center of $M$. 
\\Let   $A_P^-,$ be the set  of $P$-antidominant   elements in $A_M$. More precisely, if  $\Sigma(P)$ is  the
set of roots of $A_M $ in the Lie algebra of $P$, and  $\Delta (P)$ is the
set of simple roots, one has:
$$A_P^-\>= \{a \in  A_M \vert \vert \alpha(a )\vert _{\F} \leq 1, \alpha\in \Delta(P)\>\>\}.$$

One defines also for $\ep>0$: 
$$ A_P^-(\ep)= \{a \in  A_M \vert \vert \alpha(a)\vert _{\F}  \leq \ep, \>\> \alpha
\in \Delta (P)
\}.$$

%%%%%%%%%%%%%%%%%%%%%%%%%%%%%%%%

\noindent{\bf 2.2. Some functions on $H\bb G$.}\ste
Let $\tilde{A}_0$ be a $\si$-stable maximal split torus   of $G$, which contains a maximal $\si$-split torus $A_0$ of $G$.
Let $P_0$ be a minimal parabolic subgroup of $G$ which contains $\tilde{A}_0$. Let $K_0$ be the fixator of a special point in the apartment of $\tilde{A}_0$ in the Bruhat-Tits building of $G$. We fix an algebraic embedding 
\beq\label{to}\tau:G\to GL_n(\F).\eeq
We may and we will assume that  $\tau(K_0)\subset GL_n(\O)$ where $\O$ is the ring of integers of  $\F$ (\cite{W} I.1)). For  $g\in G$, we write:
$$\tau (g)=(a_{i,j})_{i,j=1,..,n},\>\tau (g^{-1})=(b_{i,j})_{i,j=1,..,n}.$$
We set
\beq \label{normewald}\Vert g\Vert=sup_{i,j}sup(\vert a_{i,j}\vert_\F,\vert
b_{i,j}\vert_\F).\eeq 
We have (cf. \cite{W} I.1) : 
\ber \label{nw} $\Vert g\Vert \ge 1$ for $ g\in G,\quad \Vert
g_1g_2\Vert\le
\Vert g_1\Vert\Vert g_2\Vert$ for   $g_1,\>g_2\in G$ 
and\\
$\Vert k_1gk_2\Vert =\Vert g\Vert$ for  $k_1,k_2\in K_0,\>g\in G$.\eer

We denote by 
$(\varepsilon_{M_0},\C)$ the trivial representation of the centralizer 
$M_0 $ of $\tilde{A}_0$ in $G$ and  $(\pi_0,V_0)=(i^G_{P_0}\varepsilon_{M_0},i^G_{P_0}\C)$ the normalized induced representation. Let 
$e_0$ be the unique element of
$V_0$ invariant by 
$K_0$ and such that  $e_0(1)=1$.\ste We remark that  the contragredient representation $(\check \pi_0,\check V_0)$ is isomorphic to
$(\pi_0,V_0)$. For 
$g\in G$, we set :
$$\Xi_G(g)=\langle\pi_0(g)e_0,e_0\rangle.$$
The function  $\Xi_G$ is biinvariant by $K_0$.\ste
\ber We will say that two functions $f_1$ and  $f_2$ defined on a set  $E$ with values in
$\R+$ are equivalent on a subset  $E'$ of $E$ (we  write $f_1(x)\asymp
f_2(x),\>x\in E'$), if there exist $C,C'> 0$ such that:
$$C'f_2(x)\leq f_1(x)\leq Cf_2(x),\>x\in E'.$$\eer
We recall (cf.  \cite{W},  Lemma II.1.2):
\ber\label{ldud} There exist $d\in \N$ and  for all $g_1,\>g_2\in G$, a constant  $c>0$ such that
$$\Xi_G(g_1 g g_2)\leq c \Xi_G(g)(1+ log\Vert g\Vert)^d,  g\in G.$$
 \eer
We set :
\beq\label{normeh}\Vert Hg\Vert:=\Vert \sigma(g^{-1})g)\Vert,\>g\in G.\eeq
 For a compact subset $\Omega'$ of
$G$, we deduce from (\ref{nw}):
\beq \label{ew} \Vert Hg\omega \Vert \asymp \Vert Hg \Vert,\>\omega\in
\Omega',\>g\in G. \eeq 
Let us define the functions $\Theta_G$ and  $N_d,\>d\in
\Z$ by
\beq \Theta_G(Hg)=(\Xi_G(\sigma(g^{-1})g))^{1/2},\>g\in G.\eeq
and
\beq \label{nd}N_d(Hg)=(1+log\Vert Hg\Vert )^d,\>g\in G.\eeq
(\ref{ew} ) implies (with $N=N_1$):
\beq \label{equivN} N(Hg\omega) \asymp N(Hg), g\in G, \omega \in \Omega'.Ê\eeq
The next assertion follows  from the definitions and (\ref{ldud}).
 \ber \label{lem3i} There exists $d\in \N$, and for all $g_1\in G$ there exists $c>0$ such that:
$$\Theta_G(Hgg_1)\leq c\Theta_G(g) N_d (Hg), g \in G.$$ \eer 
It follows from the Cartan decomposition for $H\backslash G$ (cf. \cite{BeOh} Theorem 1.1)  that there  exist a  compact subset $\Omega$ of $G$  and a finite set $ \mathcal{P}$ of  minimal $\si$-parabolic subgroups of $G$ such that:
\beq \label{cartan} H\backslash G= \cup_{P \in \mathcal{P}} HA^{-}_P \Omega. \eeq

\ber \label{lem1}Let $P=MU $ be a minimal $\si$-parabolic subgroup of $G$ and let $\Omega'$ be a compact subset of $G$. We  choose a norm on $\a_M$. By (\cite{L}, Lemma 7 and Proposition 6), we have:  \\
 (i) There exist $c, c',C, C'>0$ such that:
$$Ce^{c\Vert H_{M}(a)\Vert}\leq \Vert 
Ha\omega\Vert \leq
C'e^{c'\Vert H_{M}(a)\Vert},\space \omega\in
\Omega',\>a\in A_P^-,$$
\\ (ii) $$ N(Ha\omega) \asymp (1+ \Vert H_{M} (a)\Vert ), a \in A_M, \omega \in \Omega.$$
 \eer
  
\ber \label{Theta}(iii) The function $\Theta_G$ is right invariant by $K_0\cap\sigma(K_0).$
\ste
 (iv) There exist $C,C'>0$ and $d,d'\in \N$ such that for $g=
a\omega$ with $\omega\in \Omega',\>a\in A^-_P$, one has  $$C\delta_{P}^{1/2}(a)N_{-d}(Ha)\leq \Theta_G(Hg)\leq C' \delta_{P}^{1/2}(a)N_{d'}(Ha)\>.$$
\eer
\begin{Lemma} \label{conint} Let $dx$ be a non zero $G$-invariant measure on $H \bb G$. 
There exists $d\in \N$ such that:
$$\int_{H\bb G} \Theta^2_G (x)N_{-d} (x) dx < \infty.$$
\end{Lemma}
\dem
Let $P=MU\in \mm{P}$ and $\Omega $ as in (\ref{cartan}).  From (\ref{Theta}) one deduces that there exist $C'>0$ and $d'\in \N$ such that :
$$\Theta_G(H a\omega) \leq  C' \delta_{P}^{1/2}(a)N_{d'}(Ha), a \in A^-_P, \omega \in \Omega.$$
We can choose $\Omega$ large enough in order to have
$$A^-_P\Omega \subset \LL^-_P\Omega,$$
where $\LL^-_P$ is the set of  $P$-antidominant elements in $\LL(A_M)$. 
It follows from \cite{KT2},  Proposition 2.6, that  
\ber \label{kt} There exist constants $C_1,  C_2>0$ such that:$$  C_1 \delta_P^{-1}(\l) \leq {\rm vol} (H\bb H\l \Omega )\leq C_2\delta_P^{-1}(\l), \l \in \LL^-_P, $$
where $ {\rm vol}(H\bb H\l \Omega )$ is the volume of the subset $H\bb H\l \Omega $ of $H\bb G$. \eer
\\From (\ref{lem1}) (ii) one deduces that for $d''\in \N$ large enough:
$$\sum_{\l \in \LL(A_M)} N_{-d''} (H\l)<\infty.$$
This implies easily the Lemma.\qed

\section{Tempered functions on $H\bb G$}
\setcounter{equation}{0}
\noindent {\bf 3.1. On the Cartan decomposition and lattices.}\ste
Let $P=MU$ be a  $\si$-parabolic subgroup of $G$.
Let $\Sigma(P)$ be the set of roots of $A_M$ in the Lie algebra of $U$ and let $\Delta(P)$ be the set of simple roots. It will be viewed as a subset  of $\a_M'$.
Let us denote  by $ ^+{\overline{\a}_P}'$ (resp., $^+\a_P'$)   the set  of $ \chi \in \a'_M$ of the form:
$$\chi= \sum _{\alpha \in \Delta(P)} x_\alpha \alpha,$$
 where $x_\alpha \geq 0$  (resp., $x_\aa>0$) for all $\alpha \in \Delta(P)$.\\
 Let us assume that $P$ is 
 a minimal $\si$-parabolic subgroup of $G$.  If $Q=LV $ is a $\si$-parabolic subgroup of $G$ such that $P\subset Q$, let $\Delta^L$ be the set of elements of $\Delta:= \Delta(P)$ which are roots of $A_M$ in the Lie algebra of $L$. We remark  that   $A^-_Q $ is equal to  the intersection  $A_{L}  \cap A^-_P$. For $\ep>0$, we define 
$$A_P^- (Q,\ep ):=\{ a \in A^-_P\vert \vert \aa (a)\vert_\F \geq \ep, \aa \in \Delta^L \>\> and \>\> \vert \aa (a)\vert_\F<\ep, \aa\in  \Delta  \setminus \Delta^L \}.$$
Let $\mathcal{P}(P)$ be the set of $\si$-parabolic subgroups of $G$ which contain $P$.  
For $\varepsilon>0$, one has a partition of $A_P^-$:
\beq \label{partition} A^-_P=\cup_{Q\in \mathcal{P}(P)}A_P^- (Q,\ep ). \eeq 
Moreover for any $Q \in \mathcal{P}(P)$  there exists a compact subset $\omega_{\ep, Q}$  of $A_M$ such that:
\beq A_P^- (Q,\ep ) \subset  A^-_Q \omega_{\ep, Q},\\
\eeq
and further, introducing  $\Lambda^-_Q$ the set of  the $Q$-antidominant elements of $ \Lambda(A_L)$, there exists  a compact set $\omega'_{\ep, Q}$ of $A_M$ such that :
\beq \label{a+tq}A_P^- (Q,\ep ) \subset  \Lambda^-_Q \omega'_{\ep, Q}. \eeq
One uses (\ref{Mweight}) and one introduces a  multiple by $k\in\R^{+*}$ of the coweight lattice.\\
Let $\delta_\alpha\in \a_M, \alpha \in \Delta(P)$, the fundamental coweights.\\
Then $\Lambda(A_{L})$ contains a sublattice $\Lambda'_L$ of finite index in $\Lambda(A_{L})$, which is generated by $\delta'_\aa:= -k\delta_\alpha \in \Lambda^-_Q$, $\aa\in \Delta\setminus \Delta^L$ and 
by  $\Lambda(A_G)$. Let $\omega_1, \dots, \omega_p$ be a basis of $ \Lambda(A_G)$ and $ \omega'_1, \dots, \omega'_p$ be the dual basis in $\a'_G$.
 Let $\Lambda'^-_Q $ be the semigroup generated by the $ \delta'_\alpha,\aa\in \Delta\setminus\Delta^L$ and $\LL(A_G)$, i.e. :
$$\Lambda'^-_Q= \{\prod_{\alpha \in  \Delta\setminus \Delta^L}(\delta'_\alpha )^{n_\aa}\vert n_\aa \in \N\}\LL(A_G).$$Then we will see that there exists a finite set  $F_Q$ in $\Lambda'_L$ such that
\beq \label{fq}\Lambda^-_Q \subset \Lambda'^-_Q F_Q.\eeq
In fact if $\lambda \in \Lambda^-_Q$, for each $\alpha\in  \Delta\setminus \Delta^L$ (resp., $j=1, \dots, p$), one defines $n_\aa$ (resp., $n_j$)  the largest   integer such that $\langle \lambda, \aa\rangle $ (resp., $\langle  \lambda, \omega'_j\rangle $)  is less  than or equal to $-kn_\aa$ (resp.,  $n_j$).
Then  $\lambda'= \prod_{\alpha \in \Delta\setminus \Delta^L}(\delta'_\alpha )^{n_\aa}\prod_{j=1, \dots p}  \omega_j^{n_j}$ is in $ \Lambda'^-_Q$. Moreover $\l (\l' )^{-1}$ lies in a bounded subset of $ \LL(A_Q)$, as $ \lambda $ varies in $ \Lambda^-_Q$, hence  it lies in a finite set $F_Q$.\\
Summarizing, one sees that there exists a compact subset $\omega''_{\ep,Q}$ of $A_M$ such that :
\beq  \label{a+tq'}A_P^- (Q,\ep ) \subset  \Lambda'^-_Q \omega''_{\ep, Q}. \eeq
\noindent{\bf 3.2.  $\mathcal{A}(H\bb G)$, $\mathcal{A}_{temp }(H\bb G)$,  $\mathcal{A}_{2}(H\bb G)$.}\ste
The proof of the following Lemma is analogous to the proof of \cite{D1}, Lemma 3.
\begin{Lemma} \label{avf}
Let $f$ be a function on $H\backslash G$ which is right invariant by a compact open subgroup. 
The following conditions are equivalent:
\\ (i) The $G$-module $V_f$, generated by the right translates  $\rho(g)f, g\in G$,  is admissible.
\\
(ii) There exist an  admissible representation $(\pi, V)$ of $G$, an element $v$ of $ V$  and  an $H$-fixed linear form $\xi$ on $V$ 
such that $ f=c_{\xi,v}$ where:
$$c_{\xi,v}(Hg)=  \langle \xi, \pi(g)v\rangle, g\in G.$$
(iii) The function $f$ is $ZB(G)$-finite, where $ZB(G)$ is the Bernstein's center of $G$.
\end{Lemma}
We denote by $\mathcal{A}(H\backslash G)$ the vector space of such functions.
An element of this space is $A_G$-finite, hence there exists a finite set $
{\rm Exp}(f)$ of smooth characters of $A_G$ such that $$f= \sum_{\chi\in 
{\rm Exp}(f)} f_\chi,$$
where the $f_\chi$ are non zero and satisfy for some $n\in \N^*$:
$$(\rho(a ) -\chi(a))^n f_\chi=0, a\in A_G.$$
The elements of $
{\rm Exp}(f)$ are called the exponents of $f$.\\
Let $(\pi, V)$ be a smooth representation of $G$ of finite length. Then it is a finite direct sum of generalized eigenspaces under  $A_G$.    If $ \nu $ is a character of $A_{G}$, let us denote by 
$V(\nu)$ the corresponding generalized eigenspace of $V$ and by $\xi(\nu)$ the restriction to $V(\nu)$  of any element $\xi$ of $V'$, which can be extended to an element of $V'$, denoted also $\xi(\nu)$, which is   zero on the other generalized eigenspaces. If $Ê\xi\in V'^H$, $\mathrm{Exp}(\xi)$ will denote the subset of $\nu$ such that $\xi(\nu)$ is non zero.  The elements of $\mathrm{Exp}(\xi)$ are called the $A_{G}$-exponents or exponents of $\xi$. 

For any $\si$-parabolic subgroup $P$, the constant term $f_P$ of $f$ along $P$  has been defined in \cite{L}, Proposition 2.
For an $H$-invariant linear form $\xi$ on $V $,  $j^*_P(\xi)$ has been defined in \cite{L}, Theorem 1. It is an $M\cap H$-invariant linear form on the normalized Jacquet module $j_P(V)$.  One denotes by $j_P$ the canonical projection from $V$ to $j_P(V)$. 
If $f= c_{\xi,v}$, one has the equality:
\beq \label{fpc} f_P= c_{j^*_P(\xi), j_P(v)}. \eeq 
Let us recall a property of the constant term (cf. \cite{D2} Proposition 3.7), in which one has to change  right  $H$-invariance to left invariance by changing $g\mapsto f(g)$ into $g \mapsto f(g^{-1}).$
\ber \label{ct=}Let $P=MU$ be a minimal $\si$-parabolic subgroup of $G$ and let  $Q=LV$ be a $\si$-parabolic subgroup of $G$ which contains $P$. Let $K$ be an open compact subgroup of $G$. Then there exists $\ep>0$ such that, for any right $K$-invariant element $f$ of  $\mm{A} (H\bb G)$, one has
$$f (a) = \delta_Q ^{1/2}(a) f_Q(a), a \in A_M^{-}(Q<\ep ), $$
where $A_M^{-}(Q<\ep ):= \{a\in A^-_P \vert  \vert \aa(a)\vert_\F < \ep,  \aa \in \Delta(P)\setminus \Delta^L(P)\}$.
\eer
One defines $$f_P^{ind}(g) := (\rho(g)f)_P, g \in G.$$
As the Jacquet module of an admissible representation  is admissible, one deduces from (\ref{fpc}) that the constant term $f_P$ is an element of $\mathcal{A}(M\cap H\backslash M)$. 
The union $
{\rm Exp}_P(f)$ of the set of exponents of $f_P^{ind}(g), g \in G$  is finite, as the Jacquet module of the $G$-module  generated by $f$ is of finite length. This set  is called the set of exponents of $f$ along $P$. If  $\xi$ is an $H$-fixed linear form on a smooth $G$-module of finite length, one defines similarly $
{\rm Exp}_P(\xi)= 
{\rm Exp}(j^*_P(\xi))$.\ste

One says that an element $f$ of $\mathcal{A}(H\backslash G)$  is tempered (resp., square integrable) if for every $\si$-parabolic subgroup $P$ of $G$, the real part of the elements of  $
{\rm Exp}_{P }(f)$ 
are contained in $^+{\overline{\a}_P}'$  (resp.,  $^+\a'_P$).

 We denote by $\mathcal{A}_{temp}(H\backslash G)$ (resp., $\mm{A}_2(H \bb G)$) the  subspace of tempered elements (resp., square integrable) of $\mathcal{A}(H\backslash G)$. Obviously one has:
\ber \label{a2atempg}The spaces $\mm{A}_2(H\bb G) \subset \mm{A}_{temp}(H\bb G)$ are $G$-invariant subspaces of $\mm{A}(H\bb G)$. \eer

Moreover, from  \cite{KT2}, Theorem 4.7, one deduces:
 \ber  If $A_G=\{1\}$, an element $f$ of $\mm{A}(H\bb G)$ is element of $\mm{A}_2(H\bb G)$ if and only if  it is an element of $L^2(H\bb G)$. \eer
 
Let $V$ be a smooth  $G$-module of finite length. Similarly, one says that an $H$-fixed linear form $\xi$ on $V$ is tempered (resp., square integrable) if for every $\si$-parabolic subgroup $P$ of $G$ the real part of the elements of  $
{\rm Exp}_{P }(\xi)$ 
are contained in $^+{\overline{\a}_P}'$  (resp., $^+\a'_P$).  
We denote by $V_{temp}'^H$ (resp., $V_{2}'^H$) the set of tempered (resp., square-integrable) $H$-invariant  linear forms on $V$.

\begin{Lemma}\label{tempcrit}
The following conditions are equivalent:
\\(i) The function $f$ is an element of $\mathcal{A}_{temp}(H\backslash G)$ (resp., $\mm{A}_2(H \bb G)$). \\
(ii) There exist an admissible representation $(\pi, V)$ of $G$, an element $v$ of $ V$  and an element  $\xi$ in  $V_{temp}'^H$ (resp., $V_{2}'^H$)  
such that:
$$ f(Hg)= \langle\xi, \pi(g)v\rangle, g\in G.$$
\end{Lemma}
\dem One uses Lemma \ref{avf} and (\ref{fpc}). \qed 
\begin{Definition} Let $f \in \mathcal{A}_{temp}(H\backslash G)$ and let $P$ be a $\si$-parabolic subgroup of $G$. Let $
{\rm Exp} _P^{\rm w}(f) $ (resp., $
{\rm Exp}_P^+(f)$) be the set of elements $\chi$  of $
{\rm Exp}_P(f)$ such that $\rm{Re}\chi=0$ (resp., is different from zero).  The weak constant term $f_P^{\rm w}$ of $f$ along $P$  is the sum of the $(f_P)_\chi $ where $ \chi $ varies in $
{\rm Exp}_P^{\rm w} (f)$. We set $f_P^+= f_P - f_P^{\rm w}$ and $f^{{\rm w},ind}_P(g)=(\rho(g)f)_P^{\rm w}$ for $g\in G$. 
\end{Definition} 
\begin{Lemma}  \label{fpqr}With the notation of the definition, let $P=MU$, $Q=LV$ be two $\si$-parabolic subgroups of $G$ such that $P\subset Q$. Let $R= P\cap L$.
Then one has:\\
(i) $f^{\rm w}_Q  \in  \mathcal{A}_{temp}(L \cap H\backslash L)$.\\
(ii)  $$f^{\rm w}_P= (f^{\rm w}_Q)^{\rm w}_R.$$
  \end{Lemma}
\dem 
(i) From the definition of $f^{\rm w}_Q$ and the fact that $f_Q \in \mm{A}(L \cap H \bb L)$, one sees that $f_Q^{\rm w} $ is also an  element of $ \mm{A}(L \cap H \bb L)$. The set of exponents $
{\rm Exp}_R(f_Q)$ is the disjoint union of $
{\rm Exp}_R (f^{\rm w}_Q ) $ and  $
{\rm Exp}_R (f^+_Q ) $.  
From the transitivity of the constant term (cf. \cite{L}, Corollary 1 of Theorem 3), one has  $
{\rm Exp}_R (f_Q) \subset 
{\rm Exp}_P(f).$  Hence if $\chi \in 
{\rm Exp}_R (f^{\rm w}_Q ) $, one has $\rm{Re} (\chi) \in ^+{\overline{\a}_P}'$ and $\rm{Re} (\chi)$ restricted to $\a_L$ is equal to zero. This implies $ \rm{Re} (\chi) \in ^+{\overline{\a}_R}'$. One deduces (i).
\\ Let us prove (ii). 
We have $$f_Q= f^{\rm w}_Q + f_Q^+.$$ Then by the transitivity of the constant term, one has:$$ f_P= (f^{\rm w}_Q)_R + (f_Q^+)_R= (f^{\rm w}_Q)^{\rm w}_R+ (f^{\rm w}_Q)^+_R + (f_Q^+)_R. $$
Looking to exponents, one concludes that  $$ f^{\rm w}_P = (f^{\rm w}_Q)^{\rm w}_R, \quad f^+_P =  (f^{\rm w}_Q)^+_R + (f_Q^+)_R.$$ \qed
 \newpage
  
\noindent{\bf 3.3. Families of type I of tempered functions.}
\begin{Definition}\label{defi}
Let $X(\C) $  be a complex algebraic torus. We denote by $B$ the algebra of polynomial functions on $X(\C)$. We denote by   $X$  the maximal compact subgroup of  $X(\C)$. 
A  family $(F_x)$, parametrized by $X$,   of elements of   $ \mathcal{A}_{temp}(H\backslash G)$   is called  a family of type I of tempered functions on $H\bb G$ if:\\
a) There exists a compact open subgroup $J$ of $G$ such that  for all $x \in X$, $F_x$ is right invariant by $J$. \\b) For all $g\in G$, the map $x \mapsto F_x(Hg)$ is $C^\infty$ on $X$.\\
c) For all constant coefficient differential operator  $D$ on $X$ and $x\in X$, the map $Hg\mapsto D(F_x (Hg))$ is an element of $\mm{A}(H\bb G)$. \\
d) For every $\si$-parabolic subgroup  $Q=LV$   of $G$, there exists a finite  family $\Xi_Q= \{\xi_1, \dots \xi_n\}$, with possible repetitions, of characters of $A_L$   with values in the group of invertible elements $B^\times $ of $B$,  such that:\\
 (d-i)\beq \label{prodrho}(\rho(a)-(\xi_{1}(a))(x))\dots (\rho(a)-(\xi_{n}(a))(x)).F_{x, Q}^{ind}(g)= 0, a \in A_L, g\in G, x \in X,\eeq
    (d-ii)  for $i=1, \dots, n$, the real part of  $\xi_{i}(\cdot)(x)$  is  independent of  $x\in X$ and is an element of $ ^+{\overline{\a}_Q}'$.\\ In the following, we will denote $\xi_{i,x}(a)$ instead of $(\xi_i(a))(x)$.\\
    The family $\mm{E}$ of $\Xi_Q$ will be called a set of exponents of the family $F$.

\end{Definition}
We will see later (cf. Theorem \ref{eisII'}) examples of such families related to Eisenstein integrals. 
\\
The following properties are easy consequences of the definitions.
\ber \label{Ftrans} If $F$ is a family of type I, parametrized by $X$,  of tempered functions on $H\backslash G$,  the same is true for the family $\rho(g)F$, for every $g\in G$, with the same set of  exponents.\eer

\begin{Lemma} \label{Deri} Let $F$ be a family of type I, parametrized by $X$,  of tempered functions on $H\backslash G$ and $Q=LV$ be a $\si$-parabolic subgroup. \\
(i) For all   $l\in L$,  the map $x \mapsto (F_x)_Q((H\cap L) l)$ is $ C^\infty $ on $X$.\\
(ii) If $D$ is a differential operator with constant coefficients on $X$ of degree $d$, one has:  $$[(\rho(a)-(\xi_{1,x}(a))\dots (\rho(a)-(\xi_{n,x}(a))]^{2^d}D(F_x)_Q= 0, a \in A_L, g\in G, x \in X,$$ 
(iii) One has the equality:
$$(DF_x )_Q= D (F_x)_Q, x\in X.$$
In particular $D(F_x)_Q\in \mm{A} ((H\cap L) \bb L)$\\
(iv) The family $(DF_x)$  is a family of type I with a set of  exponents given by the $2^d \Xi_Q$,where $2^d \Xi_Q$ means $\Xi_Q$ repeated $2^d$ times. 

\end{Lemma}
\dem 
(i) Using translations, it is enough to prove   that $(F_x)_Q (a), a\in A_L$ is $C^\infty$ on $X$. Let $P=MU$ be a $\si$-parabolic subgroup contained in $Q$. By  (\ref{ct=}), there exists  $\ep>0$ such that 
$$(F_x)_Q ((H\cap L) a) =\delta_Q(a)^{-1/2}  F_x (Ha), a\in A_M^-(Q <\varepsilon).$$
Hence the assertion of the Lemma is true for $a\in A_M^- (Q <\varepsilon)$. But the relation (\ref{prodrho}) applied to $a_0$ strictly $P$-dominant  instead of $a$ implies  a linear  recursion relation for the sequence $ ((F_x )_Q (aa_0^{-p}))$ which allows to compute $(F_x)_Q$ on  $ (H\cap L) A_L$ from its values on $A_M^- (Q <\varepsilon)$ (cf. \cite{D2} proof of Proposition 3.11 for details). Then (i) follows.\\
(ii) By using induction, it suffices to prove the assertion for $d=1$. In that case, we apply $D$ to (\ref{prodrho}) and then we apply the product of operators $(\rho(a)-(\xi_{1,x}(a))\dots (\rho(a)-(\xi_{n,x}(a))$ to the equality obtained. This gives the result.\\ 
(iii) We fix $x\in X$. Let $V$  be the linear span of  the set  $\{ÊD (\rho(g) F_x)\vert g\in G\}$. As $D$ and $\rho(g)$ commute, this space is invariant by right translation by elements of $G$.  The elements of $V$ are of the form $DF'_x$  with  a family $F'$ of type I  satisfying (\ref{prodrho}).

We first  show that the map $DF'_x  \mapsto D(F'_x)_Q$ is well defined. For this it is enough to prove that if $DF'_x= 0$ then $D(F'_x)_Q=0$. From  (\ref{ct=}), with the notation of (ii), one has $(F'_y)_Q((H\cap L) a)=F'_y((H\cap L) a)\delta_Q(a)^{-1/2} $ on $  A_M^- (Q <\varepsilon)$ for some $\ep>0$. Hence by derivation $D (F'_x)_Q=0$ on 
 $A_M^- (Q <\varepsilon)$. Using recursion relations as in  (ii), one gets that $D(F'_x)_Q= 0$ on $A_L$. Then using translations, one sees that $D(F'_x)_Q= 0 $ on $L$.
 
Hence  the map $DF'_x  \mapsto D(F'_x)_Q$ is well defined on $V$.
From the properties of the constant term of the $F'_x$ (cf. \cite{D2}, Proposition 3.14), it is easily seen that this map   has the characteristic properties of the constant term map  on $V$ (cf. l.c.). This proves (iii).
\\  Then (iv) follows from (ii) applied to right  translates  of $F$ by elements of $G$ and from (iii). \qed

\begin{Lemma}\label{FwI} Let $F$ be a family of type I, parametrized by $X$,  of tempered functions on $H\backslash G$.  Let $Q= LV$ a $\si$-parabolic subgroup. Then, one has \\
(i)    The family  $(F_x)^{\rm w}_Q, x \in X$ is a family of tempered functions on $(L\cap H) \bb L$ of type I.\\
(ii) Let $D$ be a differential operator on $X$ with constant coefficients. Then, one has $D(F_x)^{\rm w}_Q=(D F_x)^{\rm w}_Q, x\in X$.
\end{Lemma}
\dem
(i)  Let $a'\in A_L$ be such that $ \vert \aa (a') \vert _\F <1$ for all $\aa$ in $ \Delta(Q)$. Let $ \Xi^{\rm w}_Q$ be the set of elements $\xi$ of $\Xi_Q$  such that $\xi_x$ is a unitary character of $A_L$ for all $x\in X$. We set $\Xi^+_Q= \Xi_Q\setminus \Xi^{\rm w}_Q$. We recall that there might be repetitions in these families.
From the theory of the resultant there exist elements $R,S$ of  $B[T]$ such that:
$$R(T) \prod_{ \xi \in \Xi^{\rm w}_Q} (T- \xi(a') ) + S(T) \prod_{ \xi' \in \Xi^+_Q} (T- \xi'(a') )= b, $$
where $$b=\prod _{ \xi \in \Xi^{\rm w}_Q, \xi' \in \Xi^+_Q} (\xi(a')-  \xi'(a')).$$
We define $$\Gamma_x= S_x(\rho (a') )\prod_{ \xi' \in \Xi^+_Q} (\rho(a') -\xi'_x(a') ).$$
where $\rho$ denotes the right regular representation on  the space of functions on $(L\cap H)\bb L$.
One sees easily that, from the definition of $R, S$,  the definition of the constant term and of $ \Xi_Q^{\rm w}$:
\beq \label{s} \Gamma_x (F_x)_Q= b(x) (F_x )^{\rm w} _Q, x \in X.\eeq 
From the properties of $a'$ and the definition of $\Xi_Q^+$, one sees that $b(x)$ does not vanish for $x \in X$ and is $C^\infty $ on $X$.
Hence $$(F_x )^{\rm w} _Q= b(x)^{-1} \Gamma_x (F_x)_Q, x \in X.$$
By Lemma \ref{Deri} (i),   for $l\in L$,  the map $x\mapsto  (F_x)_Q((H\cap L)l)$  is $C^\infty$ on $X$. \\  One has to prove that  for  a differential operator $D$ with constant coefficients on $X$ and $x \in X$, $D(F_x)^w_Q\in \mm{A} ((L\cap H) \bb L)$. First, from Lemma \ref{Deri}, second part of (iii),   $D(F_x)_Q$ is an element of $\mm{A} ((L\cap H) \bb L)$. Then our claim follows by applying $D$ to the preceeding equality.   \\
 Separating the exponents of $(F_x)_Q^{{\rm w}}$ and $(F_x)_Q^+$, one deduces from (\ref{prodrho}) that:
$$\prod_{\xi \in \Xi_Q^{{\rm w}}} (\rho(a) -\xi_x (a) )(F_x)_Q^{{\rm w}}=0, a\in A_L, x \in X.$$
Similarly, if $R$ is a $\si$-parabolic subgroup of $L$, one gets a relation like (\ref{prodrho}) for $((F_x)_Q^{{\rm w}})_R$. 
Altogether this shows that  $ (F_x )^{\rm w} _Q$ is a  family of type I of tempered functions on $(L\cap H)\bb L$. 
This proves (i).\\
From Lemma \ref{Deri} (iii), $(DF_x)_Q= D (F_x)_Q$. By (i),  $(F_x)_Q^{\rm w}$ is $C^\infty$ in $x\in X$. As this is also true for $(F_x) _Q$ 
this implies that $(F_x)_Q^+=  (F_x) _Q-(F_x)^{\rm w} _Q$ is also $C^\infty$ in $x\in X$. 
Hence $$(DF_x)_Q= D(F_x)_Q^{\rm w} + D(F_x)_Q^+.$$
But the exponents of $D(F_x)_Q^+$  are (up to multiplicities ) the exponents of $(F_x)_Q^+$ and similarly for $D(F_x)_Q^{\rm w} $ (cf.
 Lemma \ref{Deri} (ii)). From the definition of the weak constant term one deduces (ii).  \qed

%%%%%%%%%%%%%%%%%%%%%%%%%%%%%%%%%%
 \begin{Proposition} \label{fxqw}
Let $F$ be a family of type I, parametrized by $X$,  of tempered functions on $H\backslash G$.
Then, there exist $d \in \N$ and $C>0$ such that :
$$\vert F_x (Hg) \vert \leq C \Theta_G(Hg)N_d(Hg), g \in G,x\in X.$$
\end{Proposition}
\dem
By using  the Cartan decomposition (cf. (\ref{cartan})) and  a finite number of right translates of $F$, one sees, using (\ref{Ftrans}), (\ref{equivN}) and (\ref{lem3i}), that it is enough to prove for each element  $P =MU$ of $\mathcal{P}$, and each  family of type I, parametrized by $X$,  of tempered functions on  $H\backslash G$,  an inequality of this type for $a\in A^-_P$.   
Now, it follows from (\ref{ct=}) and Definition \ref{defi} a) that there exists  $\ep>0$ such that for all $Q \in \mathcal{P}(P)$ and for all $x\in X$:
\beq \label{fxae}{F_x}_{\vert A_P^- (Q, \ep)} =(\delta_Q)^{1/2}   { (F_x)_Q}_{\vert A_P^-(Q, \ep) }. \eeq 
By (\ref{partition}),(\ref{a+tq'}), we have  $A^-_P\subset \cup_{Q\in \mathcal{P}(P)} \Lambda^-_Q \omega''_{\ep, Q}$.  Using  a finite number of right translates,  (\ref{equivN}) again and the estimate (\ref{Theta}) of $\Theta_G$, it is enough  to prove that there exist $C>0$ and $d\in \N$ such that:
\beq \label{enough}\vert ({F_x})_Q(\lambda) \vert \leq C N_d (H\lambda), \lambda \in \Lambda'^-_Q. \eeq
By assumption on the real part of $\xi_{i,x}$, the eigenvalues $\xi_{i,x} (\l), \l \in \LL'^-_Q\subset \Lambda'_L$ have a modulus less or equal to 1 which does not depend on $x\in X$. 
We will see that (\ref{enough}) follows from the following Lemma applied to the lattice $\Lambda'_L$.
\begin{Lemma} \label{latt1}
Let $\Lambda$ be  a lattice with basis $\l_1, \dots, \l_q$.  If $\l= i_1 \l_1 +\dots + i_q \l_q $, we set  $\vert \l\vert =\vert  i_1 \vert  +\dots +\vert  i_q \vert  $. Denote by $\LL^+$ the set of $\l$ such that the $i_j$  are in $\N$. Let $\xi_{1,x},\dots  \xi_{n,x}, x\in X$ be  a $C^\infty$ family of characters of $\Lambda$  such that:
$$\vert \xi_{i,x} (\l_j)\vert \leq 1, x\in X.$$ 
Let $(f_x), x \in X$ be a $C^\infty$ family of functions on $\LL$ such that
$$(\rho (\l)- \xi_{1, x} (\l)) \dots(\rho (\l)- \xi_{n,x} (\l))f_x= 0, x\in X, \l \in \LL.$$
Then there exist  $C>0, d \in \N$ such that:
$$\vert f_x (\lambda) \vert \leq C (1+ \vert \l\vert )^d, x\in X, \l\in \LL. $$
\end{Lemma}
\dem
 If $i= (i_1, \dots, i_q) \in  \Z^q$ we define
$$\l^i= i_1 \lambda_1 +\dots  i_q \lambda_q.$$ Let $E_{n,q}$ be the space of maps from $\{0, \dots, n-1\}^q$ to $\C$. We fix a norm on this vector space.
To  $x\in X$, we associate the element $g_x$ of $E_{n,q} $ defined by 
$g_x (i)= f_x(\l^i), i\in \{0, \dots, n-1\}^q$.
Then (cf. \cite{D1}  before Lemma 14) there  exists a representation $\xi_{x}$ of $\LL$ on $E_{n,q} $, depending only on the family characters $\xi_{1,x}, \dots \xi_{n,x}$ of  $\LL$   and which  depends smoothly on $x\in X$, such that  
for $\lambda\in\LL$,  the eigenvalues of $\xi_x(\lambda)$ are  $\xi_{1,x}(\l), \dots\xi_{n,x}(\l)$ and  $$f_x(\l) = ((\xi_x(\l)g_x) (0, \dots, 0),\l \in \LL.$$
The eigenvalues of $ \xi_x (\l_1), \dots, \xi_x(\l_l)$ are of modulus  less than  or equal to 1. Moreover, from the smoothness of $ \xi_x$ in $x \in X$, one sees that the norms of the endomorphisms $ \xi_x(\l_i)$ are bounded by a constant independent from $x\in X$, as well as their inverse.\\
From \cite{DOp} Lemma 8.1, one sees that, for some $d'\in\N$,  the norm of $\xi_x (\l^i)$ is bounded by the product of a constant, independent of $x\in X$,  with  $ (1+\vert i_1\vert )^{d'} \dots (1+\vert i_q\vert )^{d'}$ for $i \in \Z^q$.
But the latter is bounded by 
$(1+ \vert i_1\vert  +Ê\dots +  \vert i_q\vert )^d$, with $d=d'q \in \N$.\qed 

 \noindent {\em End of the proof of the Proposition.}
 From (\ref{lem1}), we have:
\beq  N(Ha) \asymp(1+ \Vert H_M (a)\Vert ), a \in A_M.  \eeq
From   the equivalence of norms for finite dimensional vector spaces, one sees that:
\beq \label{equinorm} 1+\vert i_1\vert  +Ê\dots +  \vert i_l\vert \asymp N(H\l^i), i\in \Z^l,\eeq 
where $l$ is the rank of $\Lambda'_L$. Then the Lemma \ref{latt1} implies easily (\ref{enough}). This finishes the proof of the Proposition.\qed

We have the following Proposition.  
\begin{Proposition} \label{casscrit}Let $f\in \mm{A}(H\bb G)$. The following conditions are equivalent:\\
(i) The function $f$ is an element of $\mathcal{A}_{temp}(H\backslash G)$. \\
(ii) There exist $C>0$ and $d \in \N$ such that:
$$ \vert f(x) \vert  \leq C \Theta_G (x) N_d(x), x \in H\bb G.$$
\end{Proposition}
\dem
 (i) implies (ii) follows from the Proposition applied to $X$ reduced to one point.\\
One sees that (ii) implies (i)  is the analogue of (i) implies (ii) in \cite{W}, Proposition III.2.2.   Let us give a detailed proof.\\ Let $f$ as in (ii). Let $P$ be a $\si$-parabolic subgroup of $G$. Let us denote by $V_f$ the linear span of the set $\{\rho(g)f\vert g\in G\}$. \\
From (\ref{lem3i}) and (\ref{equivN}),  one sees that, for any element  $f'$  of $V_f$, there exist $C'>0$ and $d'\in \N$ such that:
\beq \label{vertf} \vert f'(x) \vert \leq C'  \Theta_G (x) N_{d'} (x), x \in H\bb G.\eeq 
Let $E_f:=\{(f'_P)_{\vert A_M}\vert f' \in V_f\}$. It is $A_M$-invariant and each element of $E_f$ is $A_M$-finite as $V_f\subset \mm{A}(H\bb G)$. One sees easily that the set ${\rm Exp}_P(f)$ is exactly the set of characters of $A_M$ which appear as  a subrepresentation of $E_f$. Let $\chi\in {\rm Exp}_P(f) $ and let $f'\in V_f$ such that $(f'_P)_{\vert A_M}$ transforms under $A_M$ by $\chi$.
  From (\ref{ct=}), one sees that there exists $\varepsilon >0$ such that:
$$f'(a) = \delta_P^{1/2}(a) f'_P(a), a\in A_M^-(P< \varepsilon)$$
From  this,  (\ref{vertf}) and(\ref{Theta}), one deduces  that there exists $d''\in \N$, $C''>0$ such that: 
$$ \vert \chi(a)\vert \leq C'' N_{d''}(a), a \in  A_M^-(P<\varepsilon). $$
Let $$\a_M^-(P< \ep)= \{ X\in \a_M\vert \> \vert  \aa(X) \vert < \ep, \aa \in \Delta(P)  \}$$
Writing $\vert \chi\vert = \chi_\nu$ for $\nu \in \a_M'$, one deduces from this and (\ref{lem1}) (ii), the existence of $C'''>0$ such that:
$$e^{\nu(X)}\leq C''' (1+ \Vert X\Vert )^{d''}, X \in \a_M^-(P< \ep).  $$
This implies that $\nu(X) \leq 0 $ for $X  \in   \a_M^-(P< \ep) $.  Hence  by applying  homotheties, one sees that $\nu$ is an element of $^+{\overline{\a}_P}'$. This proves that (ii) implies (i). \qed

We want to define  some kind of seminorms on the space of families of type I with good properties of comparison when looking to Levi subgroups. For this, we introduce suitable sets of $\si$-parabolic subgroups.

 Let $G_1$ be the group of $\mathbf  F$-points of an algebraic group  defined over $\mathbf F$. Let $\si_1$ be a rational involution of this group defined over $\F$ and let  $H_1$ be the $\mathbf  F$-points of the neutral component of the group of  fixed points of $\si_1$. If $G_1=(G_1)_{der}$, we choose a set $\mm{P}_{min}(G_1,\si_1)$ of minimal $\si_1$-parabolic subgroups of $G_1$ which gives a  Cartan decomposition for $H_1\bb G_1$ (cf. \ref{cartan}). In general, using (\ref{ll'}), let $\mm{P}_{min}(G_1,\si_1)$ be the set of minimal $\si_1$-parabolic subgroups of $G_1$ whose intersection with $(G_1)_{der}$ is an element of $\mm{P}_{min}((G_1)_{der},\si_1)$. Then from (\ref{carder}), this set of minimal $\si_1$-parabolic subgroups of $G_1$ leads to a Cartan decomposition for $H_1\bb G_1$. Let $\mm{P}(G_1,\si_1)$ be the set of $\si_1$-parabolic subgroups  of $G_1$ containing an element of $\mm{P}_{min}(G_1,\si_1)$ and $\mm{L}(G_1,\si_1)$ be the set of the $\si_1$-stable Levi subgroups of elements of $\mm{P}(G_1,\si_1)$. If there is no ambiguity on the involution $\si_1$, we drop it from the notation.
 
 We return to $G$ and $\si$, which induces an involution on each $\si$-stable subgroups of $G$. If $L$ is the $\si$-stable Levi subgroup of a $\si$-parabolic subgroup of $G$, we set:

   $$\mm{L}_1(L)=\mm{L}(L), \mm{L}_{i+1}(L)= \cup_{L_1 \in  \mm{L}_i(L)} \mm{L}(L_1).$$
If $L_1\in\mm{L}(L)$ is different from $L$ then $dim(A_{L_1})>dim A_L$. Hence, there exists $i_0$ such that  $\mm{L}_{i_0}(L)=\mm{L}_{j}(L)$ for $j\geq i_0$. We set $\mm{L}_\infty(L)=\mm{L}_{i_0}(L)$.
If $L\in\mm{L}_\infty(G)=\mm{L}_{p_0}(G)$ for some $p_0$, then $\mm{L}(L)\subset \mm{L}_{p_0+1}(G)=\mm{L}_\infty(G)$. 
Let us assume that $\mm{L}_i(L)\subset \mm{L}_{p_0}(G)$ for some $i\geq 1$. Then $\mm{L}_{i+1}(L)=\cup_{L_1\in\mm{L}_i(L)}\mm{L}(L_1)\subset\cup_{L_1\in\mm{L}_{p_0}(G)}\mm{L}(L_1)=\mm{L}_{p_0+1}(G)=\mm{L}_{\infty}(G)$.  This implies:
\ber\label{Linfini}    For $L\in\mm{L}_\infty (G)$, one has  $ \mm{L}_\infty(L)\subset\mm{L}_\infty(G).$\eer
For $L\in\mm{L}_\infty(G)$, we  denote by $\mm{P}_\infty(L)$ the set of $\si$-parabolic subgroups of $L$  whose $\si$-stable Levi component belongs to $\mm{L}_\infty(L)$. Similarly, we define 
 $\mm{L}_\infty(L_{der})$. 
Then \ber \label{adapt}  The map $M\mapsto M \cap L_{der} $ is a bijection between $\mm{L}_\infty (L) $ and $\mm{L}_\infty (L_{der}) $.
 
We say that $\mm{L}_{\infty}(L)$ is adapted to $L_{der}$. \eer

 %%%%%%%%%%%%%%%%%%%%%%%%%%
 %%%%%%%%%%%%%%%%%%%%%%%%%%%%
 We introduce the following ''seminorms'' on the space of families of type I. Notice that these seminorms might be infinite.   Let  $D$ be a finite set of differential operator on $X$ with constant  coefficients and $n\in\N$. If  $F$ is a family of type I parametrized by $X$, we set 
\beq\label{nuX}\nu^X(G,D,n,F)=\sup_{x\in X}\sup_{d\in D}\sup_{ g\in G } N_n(Hg)^{-1}\Theta_G(Hg)^{-1} \vert (d\cdot F_x)(g )\vert,\eeq
and  
\beq\label{muX}\mu^X(G, D,n,F)=\sup_{Q=LV\in\mm{P}_\infty(G)}\nu^X(L,D,n,(F)^{\rm w} _Q).\eeq
 
\begin{Remark}
Notice that in considering the right hand side of (\ref{muX}),  we have chosen the function $N$ on $L$ defined by:
$$N((H\cap L)l): = N(Hl), l\in L. $$  Another choice simply produces   functions  equivalent to this one, from (\ref{lem1}).\end{Remark} 

The following Proposition is the analogue of \cite{W} Lemmas VI.2.1, VI.2.3. The proof is essentially  similar but takes into account the dependence on the family $F$.  
\begin{Proposition}\label{fq+} We fix a set of exponents ${\mathcal E}$ and a compact open subgroup  $J$ of $G$. Let $Q= LV\in\mm{P}(G)$  and  $P=MU\in\mm{P}_{min}(G)$ such that $P\subset Q$. Let $\Delta= \Delta(P)$, $ \Delta^{L}= \Delta(P\cap L)Ê\subset \Delta$ and for $\delta>0$, let
$$D^{L} (\delta)= \{a\in A^-_P\vert \langle \aa, H_{M}(a)\rangle \leq - \delta \Vert H_{M} (a)\Vert, \aa \in \Delta \setminus \Delta^L \}. $$
There exists a compact subset  $C^L(\delta) $ of $A^-_P$  and  for all $n\in \N$, there  exist $\ep>0, C_n>0$ such that, for all families  $F$ of type I, parametrized by $X$,  of tempered functions on $H\backslash G$, with the given set of  exponents, and right invariant by $J$,  one has 
$$\vert (F_x)_Q ^+ ((H\cap L) a) \vert  \leq C_n \mu^X(G,1,n,F) \Theta_{L}((H\cap L) a)e^{-\ep \Vert H_{M} (a) \Vert}$$
for  $a\in D^{L} (\delta)\setminus C^L(\delta)$ and $x\in X$.
\end{Proposition}
\dem
One can assume that $Q$ is proper otherwise $(F_x)_Q^+=0$. We fix $n\in \N$.\\
Let us prove  that there exist    $t> 0, C_1>0$ and $d \in \N$ such that, if $a \in A_P^-$ satisfies $\langle \aa, H_{M} (a)\rangle \leq -t $ for $\aa \in \Delta\setminus \Delta^L$, one has for all families of the Proposition: 
\beq \label{1} \vert  (F_x)^+_Q((H\cap L)a)\vert \leq C_1  \mu^X(G, 1,n,F) \Theta_L ((L\cap H)a)N_{n+d}(Ha), x\in X.  \eeq 
By (\ref{ct=}), there exists $t>0$ (independent of $F$)  such that for $a$ satisfying the above hypothesis, one has the equality:
$$(F_x)_Q((H\cap L)a)= \delta_Q(a)^{-1/2}  F_x(a).$$
 By definition of the seminorms, one has
$$\vert (F_x)_Q((H\cap L)a) \vert \leq  \mu^X(G,1,n,F)\delta_Q(a)^{-1/2}  \Theta_G (Ha)N_n(Ha).$$
Applying the right inequality of (\ref{Theta})  to $G$ and the  left  inequality to $L$, and the equality $\delta_Q(a)^{-1}\delta_P(a)=\delta_{P\cap L}(a)$, one gets that there exist $C_2>0$ and $d\in \N$ such that:
$$\delta_Q(a)^{-1/2}  \Theta_G (Ha) N_{n} (Ha)\leq C_2 \Theta_L ((L\cap H)a) N_{n+d}(Ha).$$
One deduces from this an inequality like (\ref{1}) for $(F_x)_Q$.\\
A similar inequality for $(F_x)_Q^{\rm w} $ follows from   the definition of the seminorms.
%, we have $\vert  (F_x)^{\rm w} _Q(a)\vert \leq \Theta_L(L\cap H a) N(L\cap H a)^n \mu^X(L,1,n,F^{\rm w} _Q)\leq  \Theta_L(L\cap H a) N(L\cap H a)^n \mu^X(G,1,n,F)$. Hence, we obtain an equality like (\ref{1}) for $(F_x)^{\rm w} _Q$.  .
Hence (\ref{1}) follows by difference.\medskip
\\With the notations of the proof of lemma \ref{FwI}, let us define:
$$r_x(T):= \prod_{\xi \in \Xi^+_Q} (T -\xi_x (a')), x \in X.$$
By expanding these polynomials, one gets:
$$r_x (T) =\sum_{i=0, \dots, N}r_{i,x}T^{N-i}.$$
For all $\xi \in \Xi^+_Q$, $\vert \xi_x (a') \vert $ is independent of $x\in X$ and belongs to the interval $]0,1[$. Changing $a'$ to a suitable power, one can assume that:
\beq \label{rix} \vert r_{i, x} \vert \leq 2^{-i}N^{-1}, i=1, \dots, N-1. \eeq

Let us show the following property.
\ber \label{2}There exists $C_3>0$ such that, for all $k\in \N$ and all $a \in A_P^-$  satisfying $\langle \aa, H_{M} (a)\rangle \leq -t $ for $\aa \in \Delta\setminus \Delta^L$, one has:
$$\vert (F_x)_Q^+ (a(a')^k) \vert \leq C_3\mu^X(G,1,n,F) 2^{-k} \Theta_L((L\cap H)a) N_{n+d}(Ha).  $$
\eer 
If $N=0$, this implies that $\Xi^+_Q$ is empty, hence $(F_x)_Q^+=0$. So one can assume that $N\in \N^*$. Let 
 $$C_3= C_1 Sup\{2^k(N_{n+d} (H(a')^k) \vert k=0, \dots, N-1\}.$$
If $k<N$, (\ref{2}) follows from (\ref{1}) applied to $a(a')^k$,  from the definition of $C_3$,  from the equality 
$$\Theta_L((L\cap H)la')= \Theta_L((L\cap H)l), l\in L,$$ as $a'\in A_L$, and from the inequality
$$ N(Haa') \leq N(Ha) N(Ha'), $$
which follows easily  from the definitions (\ref{normeh} ) and (\ref{nd}). \\ Let $k\geq N$ and   let us assume that the inequality (\ref{2}) is true for $k'<k$. 
It follows from the definitions that $r_x(\rho(a')) (F_x)^+_Q =0$ for all $x\in X$, hence one gets:
$$(F_x)^+_Q(a(a')^k)=  -\sum_{i=1, \dots, N} r_{i,x} (F_x)^+_Q(a(a')^{k-i}).$$
The inequality (\ref{2}) for the left side of this equality follows from the induction hypothesis and (\ref{rix}).
\\ Let $C^L(\delta):= \{ a\in   D^L(\delta)\vert \>\Vert H_{M} (a)\Vert \leq t\delta^{-1} \} $. It  is compact.  Let $D= D^L(\delta)\setminus C^L(\delta)$. Hence one has:
$$D= \{ a\in D^L(\delta)\vert\> \Vert H_{M} (a)\Vert > t\delta^{-1} \}.$$
For $a\in D$, let $k$ be the largest integer  which is less or equal to 
$$ (\delta \Vert H_{M} (a) \Vert -t) (-\langle \aa, H_{M} (a')\rangle)^{-1},$$ 
when $\aa$ varies in $\Delta\setminus \Delta^L$. From the definition of $D$ and the choice of $a'$, $k$ is an element of $\N$.
From the definition of $D^L(\delta)$, $a (a')^{-k} $ is in $A^-_P$  and satisfies: \beq \label{aat} \langle \aa, H_{M} (a (a')^{-k})\rangle  \leq -t, \aa \in \Delta\setminus \Delta^L. \eeq  By applying (\ref{2}) to $a (a')^{-k} $ instead of $a$ and to the integer $k$, one gets :
$$(F_x)^+_Q((L\cap H)  a) \leq C_3 \mu^X(G,1,n,F) 2^{-k}\Theta_L ((L\cap H) a (a')^{-k}) N_{d+n}(Ha (a')^{-k}). $$
As it was already observed $\Theta_L ((L\cap H) a (a')^{-k})= \Theta_L((L\cap H) a)$.  From (\ref{aat}), one deduces:
$$Ê\vert \langle \aa, H_M(a'^k)\rangle\vert  \leq t +  \vert \langle, \aa, H_M(a)\rangle\vert.$$
From this and  (\ref{equivN}), one sees that there exists $C_4>0$ 
$$N_{d+n} (Ha (a')^{-k}) \leq C_4 (1+ \Vert H_{M} (a)\Vert )^{d+n}.$$
Writing that $$  (\delta \Vert H_{M} (a) \Vert -t) (-\langle\aa, H_{M} (a')\rangle)^{-1} \leq k+1,$$
for some $\aa \in \Delta\setminus \Delta^L$, one sees that there exist $r>0$ and $l\in \N$, independent of $a \in D$, such that:
$$r \Vert H_{M}(a) \Vert  \leq l+k.$$
From this it follows that for $a\in D$:
$$ (F_x)^+_Q(H a) \leq C_3C_4 \mu^X(G,1,n,F)\Theta_L ((L\cap H) a )2^{-r\Vert H_{M} (a)\Vert+l}(1+ \Vert H_{M} (a)\Vert )^{d+n}.$$
In order to finish the proof of (ii), it is enough to remark that there exist $C_5>0$ and $\ep>0$ such that for all $x>0$,
$$2^{-rx+l} (1+x)^{d+n} \leq C_5 e ^{-\ep x}.$$
\qed 

%%%%%%%%%%%%%%%%

\section{Wave packets in the Schwartz space}
\begin{Definition} \label{def1}
The Schwartz space $\mathcal{C}(H\backslash G)$ is the space of functions $f$ on $H\backslash G$, 
which are  right invariant by a compact open subgroup of $G$ and such that for any $d\in \N$, there exists a constant 
$C_d>0$ such that:
$$\vert f(x) \vert \leq C_d\Theta_G(x) (N_d(x))^{-1}, x \in H\backslash G.$$ 
The smallest constant $C_d$ is denoted  by $p_d(f)$.
It defines a seminorm on $\mathcal{C}(H\backslash G)$.
\end{Definition}
\begin{Lemma} One has
$$ \mm{A}_{2}(H\bb G)\subset \mathcal{C}(H\bb G). $$
\end{Lemma}
\dem  One proceeds as in the proof of  Proposition \ref{fxqw} with $X$ reduces to a single point. 
 One has to replace Lemma \ref{latt1} by the following property, which follows from [DOp], Corollary 8.2 (ii). \ber Let $A$ be an endomorpism of a finite dimensional normed vector space whose eigenvalues are of modulus strictly less than 1. Then for any $d$ in $\N$, there exists a constant $C>0$ such that:$$\Vert A^n\Vert \leq C (1+n)^{-d}, n\in \N.$$ \eer 
 This achieves to prove the Lemma. 
 \qed

\begin{Lemma} \label{AC}
  If $f\in  \mm{A}_{temp}(H\bb G)$ and $f'\in \mm{C}(H \bb G)$ the integral
$$\int_{H\bb G} f(x) f'(x) dx $$
converges absolutely.
\end{Lemma}
\dem
 The lemma  follows from  Proposition  \ref{casscrit} and Lemma \ref{conint}. \qed

\setcounter{equation}{0}
Let $M$ (resp., $L$)  be the $\si$-stable Levi subgroup of a $\si$-parabolic subgroup $P$ (resp., $Q$) of $G$.   Let $\tilde{A}$ (resp., $\tilde{A'}$) be a maximal split torus of $M$ (resp.,  $L$) such that the set  $A$ (resp., $A'$) of its antiinvariant elements is a maximal $\sigma$-split torus of $M$ (resp.,  $L$) .  By (\cite{CD} (4.9)), we can choose     a set of representatives  $W(Q\bb G/P)$  of $Q\bb G/P$  such that its elements satisfy 
    $w. \tilde{A}=\tilde{A'}$.
    
     Let $(Q\backslash G/P)_\si$ be the set of $(Q,P)$-double cosets in $G$ having a representative $w$ such that $w.A=A'$ and $w.\tilde{A}=\tilde{A'} $. 
   We denote by $W(L\backslash G /M)_\si$  a set of representatives of $(Q\backslash G/P)_\si$ with this property and we assume that $W(L\backslash G /M)_\si \subset W(Q\backslash G /P)$.
   
 Let $(Q\vert G \vert P)_\si $  be    the  set of  elements of $(Q\bb G / P)_\si$ having a representative  $w$ such that $w.\tilde{A}=\tilde{ A}', w.A= A'$ and $A_L \subset w.A_M$. Let $W (L\vert G\vert M)_{\si}$ be    a set of representatives of  $(Q\vert G \vert P)_\si $ with these properties and we assume $W (L\vert G\vert M)_{\si}\subset W(L\backslash G /M)_\si$. We want to identify $W (L\vert G\vert M)_{\si}$ with a set independent of choices.

First we prove the following property. 
\ber  \label{ssWLG} Let  $s, s' \in W(L\vert G\vert M)_\si$  such that  
  $$(s.\chi)_{\vert A_L}= (s'.\chi)_{\vert A_L}, \chi \in X(M)_{\si,u}. $$
Then  $s=s'$.
\eer

 As conjugacy by  $s$ defines an isomorphism from $A$ to $A'$, it determines a linear isomorphism   $s: \a \to \a'$. One has a similar map for $s'$. 
As  $A_L$ is contained in $s.A_M$ and in $s'.A_M $, one has   $s^{-1}\a_L  \subset 
\a_M$ and $s'^{-1}\a_L  \subset 
\a_M.$
The condition  
 (\ref{ssWLG})  implies that for all  $\lambda$  in $\a'_M$, $(s \lambda) _{\vert \a_L} =(s' \lambda) _{\vert \a_L}$. Evaluating in $X \in \a_L$, one deduces   $$ \langle\lambda, s^{-1} X\rangle=   \langle\lambda, s'^{-1} X\rangle, X \in \a_L.$$
  This implies, by varying  $\lambda$ in $ \a'_M$,
 the equality $$ s^{-1} X = s'^{-1} X, X \in \a_L.$$
 In other words, one has: 
$$s's^{-1} X= X,  X \in \a_L. $$
This  gives $s's^{-1} \in L$ and  $s, s'$ are representatives of the same  $(Q,P)$-double coset. Hence  one has  $s=s'$. This achieves to prove  (\ref{ssWLG}).\medskip

Let us remark that one has the following immediate corollary of the proof of 
 (\ref{ssWLG}).
 \ber \label{ssmu} Let  $s, s' $ be distinct elements of $ W(L\vert G\vert M)_\si$ and  $\mu$, $\mu'$ two characters  of  $A_L$. For $\chi$ in an open subset of $X(M)_{\si,u}$, one has :
 $$\mu (s.\chi)_{\vert A_L}\not= \mu' (s'.\chi)_{\vert A_L}.$$
\eer
 
 Now, we will identify $W(L\vert G\vert M)_\si$ with a set which does not depends on choices.\\
   Let $N(A,A')_\si$ be the set of $g \in G$ such that $g.A= A'$, $g.\tilde{A} =\tilde{A'}$. 
Let $M_0$ (resp., $L_0$)   be the centralizer of $\tilde{A}$ (resp.,  $\tilde{A'}$). Let $W(A,A')_\si$ be the quotient $N(A, A')_\si / M_0$ which is identified with 
 $L_0  \bb N(A, A')_\si / M_0$. It appears as  a set of isomorphisms between $\tilde{A}$ and $\tilde{A'}$. \\
 Let us prove the following result:
 
\ber \label{WWL}  Let  $(W^L)_\si $ be the subgroup of the Weyl group $W^L$  of  $L$ for $\tilde{A'}$ whose elements  preserve $A'$.   Let  
$$\overline{W} (L\vert G\vert M)_\si : =(W^L)_\si \bb \{s\in W(A,A')_\si \vert  A_L \subset s.A_M\}.$$
The natural map from   $W(L\vert G\vert M)_\si $ to  $\overline{W} (L\vert G\vert M)_\si $ is bijective.\eer

Let us prove  that this map is surjective. Let  $s\in \overline{W} (L\vert G\vert M)_\si $   and let  $s_1$  be a representative of $s$  in  $N(A,A')_\si $. 
In the   $(Q,P)$-double coset  $Qs_1P$,  there is an element $s'$ in  $W (L\vert G\vert M)_\si $ by definition of the latter set. Then one has two  elements $s', s_1\in N(A,A')_\si $ such that  $QsP=Qs_1P$. We want to show that $s'= l s_1 m $ for some  $l\in L$ and $m\in M$. 
Using conjugacy by an element of     $N(A,A')_\si $, one  can  reduce  to the case  $\tilde{A}=\tilde{A'}$. 
  But, by the Bruhat decomposition,    $s'$ and $s_1$  represent elements of the Weyl group  which have the  same  $(W^L,W^M) $-double coset. This proves the existence of $l$ and $m$.  As  $s_1. M \subset L$  one can omit  the  $m$  and write  $s'= ls_1$ for some $l\in L$. From the properties of $s', s_1$, one deduces that $l$ normalize $\tilde{A'}$ and $A'$.  Hence the image of $s'$ by our map  is $s$. Hence this map is surjective.
\\ Let us prove the injectivity. If $s, s' \in W(L\vert G\vert M)_\si $ have the same image by our map, they satisfy the condition  (\ref{ssWLG}). Hence they are equal. This achieves to prove   (\ref{WWL}).\medskip

We recall that  $X(G)_{\si,u}$ has been identified with a subgroup of $X(M)_{\si,u}$. Let $X(M)_\si^G$ (resp.,   $X(M)_{\si,u}^G$) be the neutral connected component of the group  of elements  $\chi$ of $X(M)_\si$ (resp.,  $ X(M)_{\si,u}$) whose restriction to $A_G$ is trivial. The group  $ X(M)^G_{\si, u}$ is the maximal compact subgroup of the algebraic complex torus $X(M)^G_{\si}$ and its  Lie algebra is equal to $(i\a_M^G)'$. Hence one has $X(M)_{\si,u}=X(G)_{\si,u}X(M)_{\si,u}^G$ and  $X(M)_{\si,u}^G\cap X(G)_{\si,u}$ is finite. 

Let $X$ be the maximal compact subgroup   of a complex algebraic  torus $X(\C)$. We assume  that $X(\C)$ is a finite covering of $X(M)_{\si}$ i.e.  there exists a surjective morphism of algebraic groups $p:X(\C)\to X(M)_{\si}$ whose kernel is finite. Let  $X_G$ be the neutral connected component of  $p^{-1}(X(G)_{\si,u})$ and  $X'$ be  the neutral connected component of    $p^{-1}(X(M)^G_{\si,u})$.  Then  $X_G$ (resp.,  $X'$) is  the  maximal compact subgroup of the  complex algebraic torus
equal to the connected component of $p^{-1} (X(G) _{ \si})$ (resp.,  $p^{-1} (X (M)^G_\si$). One has $X=X_GX'$ and  $X_G\cap X'$ is finite. 
For $x\in X$, we set $\chi_x=p(x)\in  X(M)_{\si,u}$.
\begin{Definition}\label{typeII'} Let $F$ be  a family of type I of tempered  functions on $H\bb G$  parametrized by $X$. \\
(1) The family  $F$ is 
called an $M$-family of type I'  if 

(i) there exists a unitary character $\mu_G$ of $A_G$ such that  
\beq \label{fchiga}F_x (Hga) =\mu_G(a)  \chi_x (a) F_x (Hg), a\in A_G, g\in G, x \in X, \eeq

(ii) Let $Q=LV$ be a $\si$-parabolic subgroup. There exists a finite set $\Xi^Q$, independent of $x\in X$, of characters of $A_L$  with $Re(\xi)\in ^+{\overline{\a}_L}' $ such that all exponents of $F_x$ along $Q$ are of the form 
$$\mu(w.\chi_x)_{\vert A_L}, \mu\in\Xi^Q, w\in W(Q\bb  G / P).$$

\noindent (2) An $M$-family of tempered functions on $H\bb G$ of type I' is said to be of type II' if for any $Q$ as above
\beq \label{fchiw} (F_x)^{{\rm w},ind}_{Q}(g)=\sum_{s\in \overline{W} (L\vert G\vert M)_\si } (F_{Q, s}(g)) _{ s.x}, x\in X, g\in G, \eeq 
where  for all $s \in \overline{W} (L\vert G\vert M)_\si $, $F_{Q,s}(g) $ is a $s.M$-family of type I' of tempered functions on $(L\cap H \bb L)$ parametrized by $s. X:= \{ (s,x)\vert x \in X\}$ with the multiplication induced by the multiplication on $X$ and with a canonical  projection on  $X(s.M)_\si $ given by $s.x:= (s,x) \to s.\chi_x$. \\
From the definition it follows that if $F$ is of type II' and $g\in G$, $\rho(g)F$ is also of type II'.
\end{Definition}
We will give examples of such families, derived from Eisenstein integrals (cf. Theorem \ref{eisII'}). Condition  (\ref{fchiw}) is motivated by this example. 
\\ Let us remark  that (\ref{ssmu})  implies the unicity of $(F_{Q,s})_{s.x}$ for $x$ in an open dense subset of  $X$ and then everywhere by continuity.  \\
 Let us prove the following assertion in which one sets $F_{Q,s}:=F_{Q,s}(e)$: 
\ber \label{FQR}  Let   $Q=LV$ be a  $\si$-parabolic subgroup of   $G$ and $s \in  \overline{W} (L\vert G\vert M)_\si$. Let  $R=SN$ be a   $\si$-parabolic subgroup of   $L$ and $s' \in \overline{W} (S \vert L\vert s. M)_\si$. Let  $Q_R= R V$. Then  $s's \in \overline{W} (S\vert G\vert M)_\si$ and 
$$\big((F_{Q,s} )_{s.x}\big)_R^{\rm w}=\sum_{s' \in \overline{W} (S \vert L\vert s. M)_\si} [(F_{Q,s})_{R, s'}]_{s's.x} $$
 with 
 $$[(F_{Q,s})_{R, s'}]_{s'sx} = (F_{Q_R,s's} )_{s's.x}.$$
 \eer 
 
 To prove this   one uses (\ref{fchiw})  for  $ Q_R$ to express directly  $F_{Q_R}^{\rm w}$, involving the second member of the equality to prove. Then, one uses  (\ref{fchiw})  for  $Q$ and $R$  and the transitivity of the weak constant term 
 (Lemma \ref{fpqr} (ii)) to compute in another way  $F_{Q_R}^{\rm w}$. Then (\ref{FQR}) follows from the identification of the terms with the same action of $A_S$ using (\ref{ssmu}). 
 \\This implies easily that:
 \ber  \label{fqs2} If $F$ is a family, parametrized by $X$, of type II' on $H\bb G$, then  $F_{Q,s}$  is a family, parametrized by $s.X$,  of type II' on $H\cap L\bb L.$ \eer 

Let $\DD^X$ be  the set of finite families of invariant differential  operators with constant coefficients on $X$. If $D\in \DD^X$ and $\psi \in C^\infty(X)$, we define:
$$q(D, \psi)= Sup \{ \vert  d\psi(\chi)\vert \vert d \in D, \chi \in X\}.$$
For $D\in\DD^X$ and  $n\in\N$, we introduce the following ''seminorms''  on families of type II'.  
\beq p^X(G, D, n,F)= \sup_{Q=LV\in\mm{P}_\infty(G)}\sup_{s\in \overline{W}(L\vert G\vert M)_\si}\nu^{s\cdot X}(L,D^s,n,F_{Q,s}),\eeq
where  $D^s$ is the family of differential operators on  $s\cdot X$ deduced from  $D$ by the action of  $s$ and $\nu$ is defined in (\ref{nuX}).  \\
If $Q=LV\in\mm{P}_\infty(G)$ then $$\nu^X(L,D,n,(F)^{\rm w}_Q)\leq \vert \overline{W}(L\vert G\vert M)_\si\vert \sup_{s\in \overline{W}(L\vert G\vert M)_\si}\nu^{s. X} (L,D^s,n,F_{Q,s}).$$ As $\mm{P}_\infty(G)$ is finite, there exists a constant $C_0>0$ such that 

\beq\label{comparemup}\mu^X(G, D,n,F)\leq C_0 p^X(G, D, n,F).\eeq

We keep the previous notation. Let $dx$ be the Haar measure of $X$ of volume 1. For  an $M$-family $F$, parametrized by $X$,  of tempered functions on $H\bb G$ of type II' and   a $C^{\infty} $ function  $\psi$ on $X$, we define 
\beq\label{SP}\mm{W}_{\psi,F} (Hg)= \int_{X}\psi(x)  F_x (Hg) dx, g \in G.\eeq

\begin{Theorem}\label{theo0}  We fix a set $\mm{E}$ of exponents and a compact open subgroup  $J$ of $G$ and let  $k\in\N$. There exist  $D, D_0\in\DD^X$   and,    for all $n\in N$,  there exists $C>0$ such that, 
 for all $M$-family $F$ of type II'  with the given set of  exponents, and right invariant by $J$, one has
\beq\label{estimII'} \sup_{g\in G}  \vert N_k (Hg)\Theta_G(Hg)^{-1} \mm{W}_{\psi,F}(Hg)\vert \leq  C p^X(G,D, n, F) q (D_0, \psi), \psi\in C^\infty (X).\eeq
\end{Theorem}
\dem

Proceeding as in the proof of Proposition \ref{fxqw},  using a finite number of right translations, one is reduced to prove a similar statement for $g \in A^-_1$, where $P_1=M_1U_1\in\mm{P}_{min}(G)$ and $A_1$ is the maximal $\si$-split torus of the center of $M_1$.  By an argument similar to (\ref{mweight}), there exists a split torus $A'_1$ of $G_{der}$, the derived group of $G$, and a finite set $F_1$ such that $A_1=A_1'A_GF_1$. Using a finite number of translations again, one is reduced to prove   the following assertion:
\ber\label{estimpo} Let $k\in \N$ and $\mm{E}$, $J$ as in the Theorem.  Then, there exist $D, D_0\in \DD^X$, and  for all $n \in \N$, there exists a constant $C>0$,  such that,  for all $M$-family $F$ of type II'  with the given set of  exponents, and right invariant by $J$, one has, for $g \in A^{'-}_1A_G$ and $ \psi \in C^\infty(X)$,
$$ \vert N_k (Hg)\Theta_G(Hg)^{-1}\mm{W}_{\psi,F}(Hg)\vert \leq  C p^X(G,D, n, F)  q(D_0, \psi).$$\eer

\noindent We first  reduce the proof of the Theorem to the case where $G$ is semisimple and then we prove it by induction on the semisimple  rank of $H\backslash G$.

\noindent{\bf Reduction to semi-simple case.} Let us assume that the Theorem is true for the derived group $G^{der}$ of $G$. \\
As  $X=X_G X'$ and $X'\cap X_G$ is a finite abelian group, there exists $C_1>0$ such that  
$$\mm{W}_{\psi,F}(Haa_1)=C_1 \int_{X_G} (\int_{X'}  \psi(x_G x ' ) F_{x_Gx' } (Haa_1) dx') dx_G, a\in A_G, a_1\in A_1^{'-}.$$
By (\ref{fchiga}), one has 
$$\mm{W}_{\psi,F}(Haa_1)=C_1\mu_G(a)\int_{X_G}x_G(a) (\int_{X' }\psi(x_G x' ) F_{x_Gx'} (Ha_1) dx') dx_G,$$
as $\chi_{x'\vert A_G}=1$.\\
By properties of the classical Fourier transform on $X_G$, for $k\in\N$, there exists $D_1\in\mm{D}^{X_G}$ such that 
$$(1+\Vert H_G(a)\Vert )^k \vert\mm{W}_{\psi,F}(Haa_1)\vert\leq C_1  \sup_{d\in D_1} \sup_{x_G\in X_G}\vert \int_{X' }d\big(\psi(x_G x' ) F_{x_G x' } (Ha_1)\big) dx'\vert,$$
and by Leibnitz formula, there exist two families $d'_1,\ldots, d'_t$ and  $d''_1,\ldots, d''_t$ in $\mm{D}^{X_G}$ such that
\beq\label{estim1}(1+\Vert H_G(a)\Vert )^k \vert \mm{W}_{\psi,F}(Haa_1)\vert  \leq C_1 \sup_{i=1,\ldots t } \sup_{x_G\in X_G} \vert \int_{X'}d_i'\cdot \psi(x_G x' ) d_i''\cdot F_{x_G x' } (Ha_1)\big) dx'\vert.\eeq

We fix  $i\in\{1,\ldots, t\}$. For $x_G\in X_G$,  we set
$$\psi'_{x_G}(x')=d_i'\cdot \psi(x_G x') ; x'\in X'$$
and $$(F'_{x_G})_{x'}(g)= d_i''\cdot F_{x_G x'}(g), g\in G_{der}.$$
For any subgroup $I$ of $G$, we set $I'=I\cap G_{der}$. We will use the notation $G'$ instead of $G_{der}$. 
As $X(M)_\si^G$ is a finite covering of $X(M')_\si^G$, $X'$ is the maximal compact subgroup of a finite covering of $X(M')_\si^G$. 
Let us prove the following assertion~:
\ber\label{restG'} The  families $(F'_{x_G})_{x'\in X'}$ are families of type $II'$ on $H'\bb G'$ with the same set of exponents $\mm{E}'$ independent of $x_G$. Moreover, they are  right invariant by $J'$.\eer

Let $Q=LV$ be a $\si$-parabolic subgroup of  $G$.
It follows from Lemma \ref{Deri} and  (\ref{FLL})  that $(F'_{x_G})_{x'}$ is of type $I$ and its exponents  along  $Q'$ are the restrictions to $A_L\cap G'$ of the exponents of $F_{x_G x' }$ along $Q$ (with different multiplicities). By definition of type $I'$, these exponents are of the form $(\mu w. (\chi_{x_G}\chi_{x'}))_{\vert A_L\cap G'}=(\mu w. \chi_{x'})_{\vert A_L\cap G'}.$ One deduces that 
\ber\label{samecontrole} $(F'_{x_G})_{x'\in X'}$ is of type $I'$  on $H'\bb G'$ and has a set of exponents along $Q'$ independent of $x_G$.  \eer

By (\ref{fPLL})  and Lemma \ref{FwI} (ii), one has 
$$((F'_{x_G})_{x'})_{Q'}^{\rm w}=\big[d''_i\cdot (F_{x_G x' })_Q^{\rm w}\big]_{\vert L\cap G'}=\sum_{s\in \overline{W}(L\vert  G \vert M)_\si} \big[d''_i.(F_{Q, s} )_{ s.x_G\; s.x'}\big]_{\vert L\cap G'}.$$ As $(F_{Q, s} )_{ s.x_G\; s.x'}$ is of type I', the same argument as before proves that  \\$\big[d''_i.(F_{Q, s} )_{ s.x_G\; s.x'}\big]_{\vert L\cap G'}$ is of type I'. So we have proved (\ref{restG'}).\medskip

We can apply our  assumption.
In the definition of the seminorms $p^X$ for $G'$, we choose the function $N$ on $G'$ defined by $N(H'g)=N(Hg)$ for $g\in G'$ (another choice produces an equivalent function).  \\
\ber\label{recu} Let  $k_1\in\N$.  There exists $D_i, D_i' \in\DD^{X'}$, and  for all $n_1\in \N$, there exists a constant $C_2>0$ such that, for all $x_G\in X_G$ and for all $a_1\in A^{'-}_1$,  one has 
 $$\vert N_{k_1} (Ha_1)\Theta_{G'}(H'a_1)^{-1}  \mm{W}_{\psi'_{x_G},F'_{x_G}}(a_1)\vert \leq C_2 p^{X'}(G',D_i, n_1,F'_{x_G}) q^{X'}(D'_i,\psi'_{x_G}).$$\eer
Let $D=\cup_{i=1,\ldots,t}\{dd''_i; d\in D_i\}$ and $D_0=\cup_{i=1,\ldots,t}\{dd'_i; d\in D'_i\}$. One has 
$$ \sup_{i=1,\ldots t } \sup_{x_G\in X_G} p^{X'}(G',D_i, n_1,F'_{x_G}) q^{X'}(D'_i,\psi'_{x_G})\leq p^X(G',D,n_1,F_{\vert G'})q^X(D_0,\psi).$$
By (\ref{estim1}) and (\ref{recu}), one deduces for all $a\in A_G $ and  $a_1\in A^{'-}_1$

\ber\label{estim2}  $$(1+\Vert H_G(a)\Vert )^k N_{k_1} (Ha_1)\Theta_{G'}(H'a_1)^{-1}\vert \mm{W}_{\psi,F}(Haa_1)\vert$$

$$
\leq  C_1C_2 p^X(G',D,n_1,F_{\vert G'})q^X(D_0,\psi).$$\eer
By (\ref{nll}),  there exist $C_L,C'_L>0$ and $ r_L, s_L$ in $\N$  such that, for $ l\in (H\cap L') \bb L'$,
\ber\label{ThetaL} $$C_L^{-1} N_{-r_L}((H\cap L)l) \Theta_L((H\cap L) l) \leq \Theta_{L'}((H\cap L')l) $$
$$ \leq C_L' N_{s_L}((H\cap L)l) \Theta_L((H\cap L)l).$$\eer 
  Let $r_0=\sup_{L\in\mm{L}_\infty(G)} r_L$ and $C_0= \sup_{L\in\mm{L}_\infty(G)}{C_L}$. Taking the inverse of the left inequality together with the fact that $\mm{L}_\infty(G)$ is adapted to $G'$ (cf.  (\ref{adapt})),  we obtain
$$p^X(G', D,n_1,F_{\vert G'})\leq C_0 p^X(G,D, n_1-r_0, F).$$

Taking $C'_G$ and $s_G$ for $L=G$ in (\ref{ThetaL}), we choose  $k_1=k+s_G$ and $n_1=r_0+n$ in (\ref{recu}). By    (\ref{estim2}), we obtain 
$$(1+\Vert H_G(a)\Vert )^k N_{k}(Ha_1)\Theta_G(Ha_1)^{-1}\vert\mm{W}_{\psi,F}(Haa_1)\vert$$
$$\leq C_1C_2 C_0 C'_G p^X(G,D,n,F)q^X(D_0,\psi).$$
By (\ref{nw})  and (\ref{lem1}), there exists $C_3>0$ such that 
$$N_k(Ha a_1)\leq C_3(1+\Vert H_G(a)\Vert )^k N_{k}(Ha_1).$$

Recall that  $\Theta_G(Haa_1)=\Theta_G(Ha_1)$ for $a\in A_G$.  Then one deduces (\ref{estimpo}) from the previous inequalities. \medskip

\noindent{\bf Semisimple case.} We prove the Theorem by induction on $dim A_1$. If  $dim A_1=0$ then $H\bb G$ is  compact and  the result   is clear. 
Let us assume  that $dim A_1>0$.  Let $\Delta_1:=\Delta(P_1)$. If $a\in A_1$, we define $s(a)= inf \{ \langle \aa,H_{M_1} (a)\rangle \vert \aa \in \Delta_1\}$. For $\Phi\subset \Delta_1$, we define $A_1^-(\Phi)$ to be the set of all $a\in A_{1}^-$ such that $\Phi=\{\aa\in\Delta_1| \langle\aa, H_{M_1}(a)\rangle=s(a)\}$. So, one has  $A_{1}^-= \cup_{\Phi\subset \Delta_1}A_1^-(\Phi)$.   Let $Q=LV\in\mathcal{P}(G)$  be such that $P_1\subset Q$ and $\Phi=\Delta_1-\Delta_1^L$. If $\Phi=\emptyset$ (which corresponds to $Q=G$) then $A_1^-(\Phi)=\emptyset$, from our hypothesis on the semisimple rank of $H\bb G$. Hence we  can assume $Q\neq G$.\medskip
 
As the set  $\{ a\in A_1^-(\Phi)| t\leq s(a)\leq 0\}$ is a compact subset of $A_1^-$,  the inequality  (\ref{estimpo}) on  this set is clear. It is enough to prove the statement (\ref{estimpo}) on $A_1^-(\Phi,<t):=\{ a\in A_1^-(\Phi)| s(a)<t\}$ for some $t<0$. 
By  (\ref{fxae}), there exists $t< 0$, which depends only on $J$,  such that, for  all $\si$-parabolic subgroup containing $P_1$ and $a\in A_{1}^-$ with  $s(a)<t$, one has
 \ber \label{FwF+} $$F_x(Ha)=\delta_Q(a)^{1/2} (F_x)_Q((H\cap L)a)$$
$$ =\delta_Q(a)^{1/2} (F_x)^{\rm w}_Q((H\cap L)a)+\delta_Q(a)^{1/2} (F_x)^+_Q((H\cap L)a).$$\eer 
We fix such $t<0$.
\\ By  (\ref{fqs2}), one has $F_Q^{\rm w} =\sum_{s\in\overline{W}(L\vert G\vert M)_\si} F_{Q,s}$ where   $F_{Q,s}$ is of type II' on $(L\cap H)\bb L$. \\  By the   induction hypothesis applied to $L^{der}\cap H\backslash L^{der}$ and the reduction to the semisimple case, the Theorem is true for $L\cap H\backslash L$. Let $k_1\in \N$.  For   $s\in \overline{W}(L\vert G\vert M)_\si$,   there exist $D_1, D'_0\in\mm{D}^X$ and for $n_1\in \N$, there exists $C'_0>0$ such that,  for $a\in A_1^-$, 
$$\vert N_{k_1}(H a) \Theta_{L}((L\cap H) a )^{-1}\vert\mm{W}_{\psi, F_{Q,s}}(a)\vert\leq C'_0 p^{s.X}(L, D_1^s,n_1, F_{Q,s}) q(D'_0,\psi).$$\\
Recall that $$ p^{s.X}(L, D^s_1, n_1, F_{Q,s})=\sup_{R=SN\in\mm{P}_\infty(L)} \sup_{s'\in \overline{W}(S\vert L\vert s.M)_\si}\nu^{s's.X}(S,D_1^{s's}, n_1, (F_{Q,s})_{R,s'}).$$
\\ As $\mm{L}_\infty(L)\subset \mm{L}_\infty (G)$  (cf. (\ref{Linfini})),  by  (\ref{FQR}) and the finiteness of $\overline{W}(L\vert G\vert M)_\si$,  there exist $D, D_0\in\mm{D}^X$, and for $n_1\in\N$ there exists  $C_0>0$ such that
$$\vert N_{k_1}(H a) \Theta_{L}((L\cap H) a )^{-1}\vert\mm{W}_{\psi, F^{\rm w}_Q}(a)\vert\leq C_0 p^X(G,D,n_1, F)  q(D_0,\psi).$$\\
By $(\ref{lem1})$ and $(\ref{Theta})$, there exist $C_2>0$ and $r\in\N$ such that 
\beq\label{TL<TG} \delta_Q^{1/2}(a)\Theta_L((L\cap H) a)\leq C_2 N_r(Ha)\Theta_G(Ha), a\in A_1^-.\eeq
Taking $k_1=k+r$ and $n_1=n$, this gives  an upper bound  like  (\ref{estimpo}) for \\
$\vert \delta_Q(a)^{1/2} \int_{X} \psi(x) (F_x)^{\rm w}_Q(a)dx \vert $ on $  A_1^-$.\\
With the notation of Proposition \ref{fq+}, there exists $\delta>0$ such that $A_1^-(\Phi, <t)\subset D^L(\delta)\subset A_{1}^-$. By  Proposition \ref{fq+}  and (\ref{TL<TG}), for  $n\in \N$, there exist  $\varepsilon>0$ and $C_n>0$   such that, for $a\in  D^L(\delta)\setminus C^L(\delta) $, one has
$$\delta_Q^{1/2}(a) |\int_{X} \psi(x) (F_x)^+_Q(a) dx |\leq $$ $$ C_n\sup_{x\in X}\vert \psi(x)\vert \mu^X(G,1,n,F)  \Theta_G(Ha) e^{-\varepsilon \Vert H_{M_1}(a)\Vert} N_r(Ha),$$
for $a\in D^L(\delta)\setminus C^L (\delta) $.\\
 As  $ \mu^X(G,1,n,F) \leq  C p^X(G,1,n,F)$ for some constant $C>0$ and for $k\in\N$, there exists $C_k$ such that $N_r(Ha)e^{-\varepsilon \Vert H_1(a)\Vert}\leq C_k N^{-k}(Ha) $, one deduces an upper bound  like  (\ref{estimpo}) for $\vert \delta_Q(a)^{1/2} \int_{X} \psi(x) (F_x)^+_Q(a)dx\vert$  for  $a\in  A_1^-(\Phi,<t)\setminus C^L(\delta)$.  Together  with the result above for $F^{\rm w}_Q$ and (\ref{FwF+}), one gets  (\ref{estimpo}) for $a\in  A_1^-(\Phi,<t)\setminus C^L(\delta)$. \\
 As $C^L(\delta)$ is compact, one gets a similar inequality for $a\in C^L(\delta) $. This achieves the proof.\qed

\begin{Theorem} \label{theo1} Let $F$ be an $M$-family, parametrized by $X$,  of tempered functions on $H\bb G$ of type II'.
Let $\psi$ be a $C^{\infty} $ function on $X$.  
\\ (i) $\mm{W}_{\psi,F}$ is an element of $\CC (H\bb G)$.
\\(ii) For each $k \in \N$, there exists a continuous semi norm $q_k$ on $C^\infty(X) $ such that (with the notation of Definition \ref{def1}):
$$p_k (\mm{W}_{\psi,F}) \leq q_k(\psi), \psi \in C^\infty(X).$$
\end{Theorem}
\dem Let $k\in\N$. We fix $D$ and $D_0$ in $\mm{D}^X$ as in the Theorem \ref{theo0}. Let $Q=LV\in\mm{P}_{\infty}(G)$,  $s\in \overline{W}(L\vert G\vert M)_\si$ and $d\in D^s$. By assumption and Lemma \ref{Deri} (iv), $d.F_{Q,s}$ is a $s.M$ family of type I on $L\cap H\bb L$. By Proposition \ref{fxqw}, there exist $n=n(Q,s,d)\in\N$ and $C_n>0$ such that  $$\sup_{x\in X} \sup_{l\in L\cap H\bb L} \Theta^{-1}_L(L\cap H l) N_{-n}(H\cap L l)\vert d.(F_{Q,s})_x(H\cap L l) \vert\leq C_n.$$
As $\mm{P}_{\infty}(G)$ is finite as well as  $D$ and $ \overline{W}(L\vert G\vert M)_\si$  for $L\in \mm{L}_\infty(G)$, there exist $n_1\in\N$ and  $C_{n_1}>0 $ such that, for $L\in\mm{L}_\infty(G)$ and $ s\in \overline{W}(L\vert G\vert M)_\si$, one has $$\nu^{s.X}(L, D^s, n_1,   F_{Q,s})\leq C_{n_1}.$$One deduces 
$$p^X(G, D, n_1,   F)\leq C_{n_1}.$$
The Theorem follows from the Theorem \ref{theo0} with $q_k=C_{n_1} q(D_0, .)$.
\qed
%%%%%%%%%%%%%%%%%%%%%%%%%%%%%%%%%%%%

 %%%%%%%%%%%%%%%%%%%%%%%%%%%%%%%%%%
  \section{Some properties of Eisenstein integrals}
 \setcounter{equation}{0}
 \noindent{\bf 5.1. Eisenstein integrals.} \\
 Let us recall some results of \cite{BD}.
Let $P=MU$ be a $\si$-parabolic subgroup of $G$, $(\delta, E)$ be a  smooth representation of finite length of $M$.
Let  $I_{\delta}$ be the space  of the induced representation  $i_{K_0 \cap P}^{K_0} \delta _{\vert_{K_0
\cap P}}$. Let $i^G_PE_\chi$ or simply $I_{\delta_\chi}$ be the space of the normalized  induced representation  $\pi_{\chi}:=i_P^G (\delta_
\chi )$, $\chi \in X(M)_\si$, where $\delta_\chi = \delta\otimes \chi $.  The  restriction of  functions from $G$ to  $K_0$ determines an isomorphism
of  $K_0$-modules between $I_{\delta_\chi}$   and $I_\delta$. One denotes by  $ \overline {\pi}_{\chi} $  the representation of  $G$  on  $I_{\delta}$  deduced from   $\pi_{\chi}$  by '' transport de structure''  by this isomorphism.\ste If  $\varphi \in I_{\delta}$  and $\chi \in X(M)_\si$, one denotes by $\varphi_{\chi}$ the element of the space $I_{\delta_\chi}$   corresponding  to  $\varphi$ by this isomorphism.\ste 
Let $B$ be the algebra of polynomial functions  on $X(M)_\si$, which is  generated by the functions $b_m, m\in M$ defined by $b_m(\chi)= \chi(m)$. 
One has:
    \ber For all 
$\varphi \in I_{\delta}$ and $g \in G$, $\chi\mapsto {\overline\pi}_{\chi}(g) \varphi$ is an element of $I_\delta$ $\otimes$ $B$.\eer 
Let  $\mm{ O}$ be the union of the open  $(P,H)$-double  cosets in 
$G$.  There exists a set of representatives,   $\overline{\mm{W}}_M^G$,     of these  open $(P,H)$-double  cosets which depends only on $M$ and not on $P$.  Moreover  for all $x \in \overline{\mm{ W}}_M^G$, $x^{-1}.P$ is a $\si$-parabolic subgroup of $G$ (cf \cite{BD} Lemma 2.4). Let $A$ be a maximal $\si$-split torus of $M$. We may (cf. \cite{BD} , beginning of section 2.4 and Lemma 2.4 ) and we will assume that for all $x \in \overline{\mm{ W}}_M^G$, $x^{-1}.A$ is a $\si$-split torus. One says that $x$ is $A$-good. Then $x^{-1}.M$ is the $\si$-stable Levi subgroup of $x^{-1}.P$ (cf \cite{CD}, Lemma 2.2). \\ One sets $J_\chi = \{ \varphi \in I_{\delta_\chi} \vert \rm{Supp} (\varphi) \subset  \mm{ O}\}$ and we define: 
 $${\mathcal V}_{\delta}:
=  \oplus _ {x \in {\overline{\mathcal W}}_M^G} ({E'_{\delta}})^{M \cap x.H }.$$ 
Let $\chi \in X(M)_\si$. To    $\eta \in
{E'_\delta}^{M \cap x.H}$,  one associates $j (P,\delta_\chi,\eta) \in J'_\chi$ defined by:
\beq \label{defj} j(P,\delta_\chi,\eta)(\varphi) = \int_{H \cap (x^{-1}.M )\bb H} \langle\varphi(xh),\eta\rangle d\dot{h},
  {\varphi \in J_\chi}.\eeq 
  Then one has (cf.  \cite{BD}, Theorem 2.8):
 \ber \label{defxi}  For $\chi$ in an open dense subset  $O$ of $X(M)_\si$, $j(P, \delta_\chi, \eta)$ extends uniquely to an $H$-invariant linear form $\xi(P, \delta_\chi, \eta)$ on $I_{\delta_\chi}$. 
There exists a  non zero polynomial $q$ on $X(M)_\si$  such that for all 
 $\varphi \in I_\delta$, the map $\chi \mapsto
q(\chi)\langle \xi(P,\delta_\chi,\eta), \varphi_\chi\rangle$,  defined on $O$,  extends to a polynomial function on $X(M)_\si$.\eer 
The Eisenstein integrals are defined, as rational functions of $\chi\in X(M)_\si$, by
\beq \label{defeis} E^G_P (\eta\otimes \varphi_\chi )(Hg) = \langle \xi(P, \delta_\chi, \eta), \pi_\chi (g) \varphi_\chi\rangle, g\in G, \varphi \in I_\delta. \eeq
Let $x\in \overline{\mm{W}}^G_M$ and $\eta\in E'^{M\cap x.H}_{\delta}$.   Then from our choice,  $x^{-1}.P$ is a $\si$-parabolic subgroup and  $x^{-1}.M$ is its $\si$-stable Levi. One can choose 1 as an element of $\overline{\mm{W}}^G_{x^{-1}.M}$ and one has $ E'^{M\cap x.H}_{\delta}=E'^{(x^{-1}.M)\cap H}_{x^{-1}\delta} \subset \mm{V}_{x^{-1}.\delta}$.  Let $\chi \in X(M)_\si$. The map $ \varphi \to \l(x^{-1})\varphi$ is a
bijective intertwining map between $i^G_P (\delta_\chi)$ and $i^G_{x^{-1}.P} ((x^{-1}\delta)_{x^{-1}\chi})$. By ''transport de structure'', one sees 
\beq \label{epexp} E^G_P (\eta\otimes \varphi_\chi)= E^G_{x^{-1}.P}(\eta \otimes (\l(x^{-1}) \varphi_\chi)), g \in G, \varphi \in I_\delta. \eeq 
%%%%%%%%%%%%%%%%%%%%%%%%%%%%%%%%%%%%
\noindent{\bf 5.2. Examples of families of type II' of tempered functions.}\\
We keep the notation of the previous subsection (cf. also after (\ref{a2atempg})).  Let $E'(\delta,H)_{2}=\oplus_{x\in \overline{\mm{W}}^G_M} {E'_2}^{M\cap x.H} $. 
One introduces the    $C$ functions  as in \cite{CD}, Proposition 6.2 (resp., Theorem 8.4).  One has
\ber \label{cteis}If   $ \phi \in   E'(\delta,H)_{2}  \otimes i^G_P E_\chi$ and $\delta$ and  $\chi$ are  unitary, $E^G_P(\phi)$ is tempered when it is defined and   one has the following  equalites  of rational functions in $\chi$:\\
Let $Q=LV$ be a $\si$-parabolic subgroup. Then

 $E^G_P(\phi)^{\rm w}_{Q}=0$ if  $W(L \vert G\vert M)_\si$  is empty.\\
Otherwise   $E^G_P(\phi)^{\rm w}_{Q}$ is equal, with the notation of l.c.,  to: 
$$ \sum_{w \in W (L \vert G\vert
M)_\si}E^{L}_{L\cap w.P}(r_{L} (C(w, Q, P, \delta_\chi) \phi)).$$Ê\eer 

\begin{Theorem} \label{eisII'} Let $\varphi\in I_\delta$ and 
 $\eta\inÊE'(\delta,H)_{2}$.
There exists a non zero polynomial $p$ on $X(M)_\si$, such that 
$\chi \mapsto p(\chi) E^G_P(\eta\otimes \varphi_\chi), \chi\in X(M)_{\si,u}$ is 
a family  $F_\chi$    of type II' of tempered functions on $H\bb G$. 
\end{Theorem} 
\dem 

 From the definition of $C$-functions and Proposition 6.2 in l.c., one sees that, for all $\si$-parabolic subgroup $Q$, there is a non zero polynomial $p_Q$ in $\chi\in X(M)_\si$ such that its product with each   term of the sum in (\ref{cteis}) is polynomial. As the set of $H$-conjugacy classes of maximal $\si$-split tori is finite, the set of $H$-conjugacy classes of minimal $\si$-parabolic subgroups is finite. Let $\mm{F}$ be a set of representatives of these $H$-conjugacy classes.

 Let $p$ be the product of the polynomials $p_Q$ when $Q$ varies in the (finite) set  $\mm{P}$ of $\si$-parabolic subgroups of  $G$ containing an element of $\mm{F}$.    We will see that this polynomial satisfies the property of the Theorem. 
One has to check  the various properties of the Definition \ref{typeII'}. 

Let  $\mu_G$ be the restriction to $A_G$ of the central character of $\delta$. The property 1) (i) is true with $\mu_G$ as  the central character of $i^G_P \delta_\chi$ restricted to $A_G$  is the  restriction to $A_G$ of the central character of $\delta_\chi$.  The exponents of $(F_\chi)_Q$ are restriction to $A_Q$ of  exponents of the Jacquet module of $i^G_P \delta_\chi$. These have the required form (see e.g. \cite{D1}, Lemma 7 (i)  which is valid without assumptions on $\delta$ and $\chi$). 
\\ It remains to check the  property c) of the Definition of families of  functions of type I (cf. Definition \ref{defi}). \\
 For this, we need to use the notation  of Definition 4 in \cite{D1}. Namely, let $\mm{O}$ be a $C^\infty$ manifold, $V$ be a vector space and for all $\nu \in \mm{O}$,  $\pi_\nu $ be an admissible  representation  of $G$ on $V$, such that the action of some maximal compact subgroup does not depend on $\nu$. For all $v \in V$ and $g \in G$, the map $\nu\to \pi_\nu (g) v$  varies in a finite dimensional vector space of vectors   fixed by some  compact open subgroup. Let us  assume that it is $C^\infty$ for all $v\in V$.  We say that $(\pi_\nu)$ is a  $C^\infty$-family of representations of $G$ in $V$.
\\Let $D$ be a $C^\infty $ vector field on $\mm{O}$. Let us define   a smooth family  $(D. \pi_\nu)$ of representations of $G$ in $V\times V$ as in \cite{D1}, Lemma 16, by:
$$(D.\pi_\nu)(g) (v_1, v_2)= (\pi_\nu(g) v_1+D (\pi_\nu (g) v_2), \pi_\nu (g)v_2), g \in G, v_1, v_2 \in V.$$
Let $ (\xi_\nu)_{\nu\in \mm{O}}$ be a family of linear forms on $V$, such that for all $\nu$, $\xi_\nu$ is $H$-fixed by the dual representation of $\pi_\nu$, and for all $v\in V$, $\nu\mapsto \langle v, \xi_\nu \rangle $ is $C^\infty $ on $\mm{O}$. 
Let  $v\in V$ and let us denote by $F_\nu$ the generalized coefficient coefficient $g\mapsto \langle \pi_\nu(g) v, \xi_\nu\rangle$. 
Then $DF_\nu$ appears as a  generalized coefficient  of $D. \pi_\nu$. More precisely let $ \tilde{\xi_\nu} = (\xi_\nu,  D\xi_\nu)\in V'\times V'$.
Then $$DF_\nu=\langle D.\pi_\nu (g) (0, v), \tilde{\xi_\nu}\rangle.$$
A simple  computation  shows that $\tilde{\xi}_\nu$ is $H$-fixed by the dual representation  of $D.\pi_\nu$.  
This formula implies that $DF_\nu $ is in $\mm{A}(H\bb G)$. Applying an induction process, one sees that it is true for any $C^\infty$ differential operators on $\mm{O}$. 
\\ This applies immediately to Eisenstein integrals and this proves that our family $F$ satisfies condition c) of Definition \ref{defi}. Hence $F$ is of type I'.  \\ 
The decomposition of $F^{{\rm w}, ind}_Q$ as in (\ref{fchiw})  follows from (\ref{cteis}). If $Q\in \mm{P}$ each term of this sum is  of type I', as it follows from the properties of $p_Q$ and the definition of $p$.  If $Q$ is $H$-conjugate to an element of $\mm{P}$, it is also true due to (\cite{D2}, Proposition 3.16 (ii)).  Hence $F$ is of type II' and the Theorem is proved.\qed

 %%%%%%%%%%%%%%%%%%%%%%%%%%%%%%%%%%%%%%%%%%%%%
\section { \label{6} Appendix: some properties of the derived group}
\setcounter{equation}{0}
 Recall that we denote by $G_{der}$ the group of $\F$-points of the derived group  $\underline{G}_{der} $ of $\underline{G}$.  If $J$ is a subgroup of $G$ we denote, unless otherwise specified,  by $J'$ the intersection of $J$ with $G_{der}$. In particular, one has $G_{der}=G'$. 
Let $Z(G)$ the group of $\F$-points of the center of $\underline{G}$. 
We recall the following facts:
\ber \label{g'ag}The group $G' \tilde{A}_G$ is cocompact in $G$ and the group $(H\cap  \tilde{A}_G) A_G$ is of finite index in  $\tilde{A}_G$.Ê\eer 

If $\tilde{A}_0$ is  maximal split torus of $G$, there exists a maximal split torus $\tilde{A}'_0$ of $G'$  such that  $\underline{\tilde{A}}_0= \underline{\tilde{A}}_G \underline{\tilde{A}}'_0$: this has been proved for at least one  $\tilde{A}_0$ in the proof of (\ref{mweight}) and the result follows from the fact that all maximal split tori of $G$ are $G$-conjugate.   It is clear that $\tilde{A}'_0$ is the maximal split torus of $\tilde{A}_0Ê\cap G'$ and one has \ber \label{aa'} The map $\tilde{A}_0Ê \to \tilde{A}'_0Ê$ is a bijection between  the set of maximal split torus of $G$ and the set of maximal split torus of $G'$.\eer 
Hence one has:
\ber \label{ag'} All maximal split tori of $G$ are $G'$-conjugate. \eer 
If $ \l\in \Lambda(\tilde{A}_0)$, let $P_\l$ be the parabolic subgroup of $G$ which contains $\tilde{A}_0$,  such that the roots of $\tilde{A}_0$ in the Lie algebra of $P_\l$ are the roots $\aa$ such that $\vert \aa (\l)\vert_\F\leq 1$. \\
If $P$ is a  parabolic subgroup of $G$ and  $\tilde{A}_0 \subset P$, there exists $ \l\in \Lambda(\tilde{A}_0)$ such that $P=P_\l$.  One has seen (cf. after (\ref{mweight}))  that the lattice $\LL(\tilde{A'_0} )\LL (A_G)$ is of finite index in $\LL(\tilde{A}_0)$. Then  a power of $\l$ is an element of this lattice, hence of the form $\l' \mu$ where $\l'\in \LL(\tilde{A}'_0) $ and $\mu \in \LL (A_G)$. One deduces  from the definitions the equality: \beq \label{pll'} P_\l = P_{\l'}.\eeq Hence  one can even choose $ \l \in \LL(\tilde{A}'_0)$. From this and \cite{BD} Equation (2.7), it follows easily that $P\cap G'$ 
is a parabolic subgroup of $G'$. Reciprocally if $P'$ is a parabolic subgroup of $G'$ then there exists   $ \l \in \LL(\tilde{A}'_0)$ such that $P'=P_\lambda\cap G'$. Looking to Lie algebras, one sees that $P_\l$ is the unique parabolic subgroup  of $G$ such that $P'=P_\lambda\cap G'$. Altogether we have shown:
\ber \label{ppÔ} The map $P \mapsto P \cap G'$ is a bijection between the sets of parabolic subgroups of $G$ and $G'$. \eer 
If $P$ and $Q$ are opposed parabolic subgroups of $G$, one can choose   a maximal split torus $\tilde{A_0}\subset P\cap Q$ and $\l \in \Lambda(\tilde{A}_0)$ such that $P=P_\l$ and $Q= P_{\l^{-1}}$. As above we can take $\l\in  \Lambda(\tilde{A}'_0)$. This implies that $P\cap G'$ and $Q\cap G'$ are opposed parabolic subgroups of $G'$. One shows similarly that if $P'$, $Q'$ are opposed parabolic subgroups of $G'$ and $P$ (resp., $Q$) is the unique parabolic subgroup of $G$ which contains $P'$ (resp., $Q'$) then $P$ and $Q$ are opposed.

 It follows easily  that the  map $P\mapsto P'= P\cap G'$  is a bijection between the sets of $\si$-parabolic subgroups of $G$ and $G'$,  and in particular between the sets of minimal $\si$-parabolic subgroups. Then  it follows: \ber \label{mm'} The  map $M \mapsto M \cap G'$ is a bijection for the sets  Levi subgroups of $\si$-parabolic subgroups of $G$ and $G'$, \eer
 which can be  specialized to   Levi subgroups of minimal $\si$-parabolic subgroups.
 The map which associates to such a Levi subgroup its unique  maximal $\si$-split torus is a bijection (cf. \cite{HW} Proposition 4.7 and Lemma 4.5).  
 Hence it follows that the correspondence  which associates to a maximal $\si$-split torus $A$ of $G$  the maximal split torus  $A'$ of its intersection with $G'$ is a bijection between the sets of maximal $\si$-split tori of $G$ and $G'$.  Then one has:
 \ber  \label{uniquesplit} The split torus $A$ is the unique maximal $\si$-split torus such that $A'\subset A$. \eer This implies, for reason of dimensions, that $\underline{A}= \underline{A'} \underline{A_G}$.  From which it follows:
 \ber \label{a1g'a} Let   $A_1$ be a maximal $\si$-split torus. If  $A'_1= g'. A'$, for some $g'\in G'$, one has $A_1= g'.A$. \eer 
 Hence it follows  from (\ref{ag'}) that: 
\ber  All the maximal $\si$-split tori of $G$ are $G'$-conjugate.\eer 
 Let $P=MU$ be a $\si$-parabolic subgroup of $G$. Recall that $M':=M\cap G'$. Let us show that  \beq \label{asubset} A_{M'} \subset A_M. \eeq  
 One has only to check that $A_{M'}$ is in the center of $M$. But the derived group $\underline{M}_{der}$ of $\underline{M}$  is contained in $\underline{G}_{der}$, hence contained in $\underline{M'}$. As    $\underline{M}$ is the almost product of  $\underline{M}_{der}$ and its center, an element of  $\underline{M}$ which commutes with  $\underline{M}_{der}$ is an element of the center. Our claim follows easily. \\
  There exists $\l\in \LL(A_M)$ such that $P=P_\l$. As in the proof of (\ref{ppÔ}), one shows that:
 \ber \label{pp'si} There exists $\l\in A_{M'} $ such that $P=P_\l$. \eer 
 Applying (\ref{ppÔ}) and (\ref{mm'})  to a $\si$-stable  Levi subgroup $L$ of a $\si$-parabolic subgroup of $G$  and to $L'$ one sees:
  \ber \label{ll'} The map $P\mapsto P\cap L'$ is a bijection between the sets of $\si$-parabolic subgroups of $L$ and $L'$. 
 The  map $M \mapsto M \cap L'$ is a bijection for the sets of  Levi subgroups of $\si$-parabolic subgroups of $L$ and $L'$, \eer
%%%%%%%%%%%%%%%%%%%%
As the derived group of $L$ is contained in $L'$, by  (\ref{g'ag})   applied to $L$ instead of $G$, one deduces that 
 \ber\label{carder}  If a finite set of minimal $\si$-parabolic subgroups  of  $L'$ leads to a Cartan decomposition for $(L'\cap H)\bb L'$, the corresponding family  of minimal $\si$-parabolic subgroup of $L$ leads to a Cartan decomposition for $(L\cap H)\bb L$. \eer
 \\
  Let us prove: 
\ber  \label{fPLL}  If $f\in \mm{A}((L\cap H) \bb L)$ then  $f_{\vert L'}\in \mm{A}((H\cap L')\bb L')$ and  for  any $\si$-parabolic subgroup $P=MU$ of $L$, one has
$ (f_P)_{\vert M\cap L'}=(f_{\vert L'})_{P\cap L'}$. Moreover, if  $f \in \mm{A}_{temp}((H\cap L)\bb L)$    then  $f_{\vert L'}\in \mm{A}_{temp} ((H\cap L')\bb L')$  and $(f_P)^{\rm{w}}_{\vert M\cap L'}=(f_{\vert L'})^{\rm{w}}_{P\cap L'}$. \eer
It follows from Lemma 2.1 in \cite{GK}, that   a   finitely generated admissible $L$-module is also an  admissible finitely generated $L_{der} $-module. Hence the same property is true for $L'$. This implies easily the former half of our first claim. \\ The space $V$ of restriction to $(L'\cap H)\bb L'$ of smooth functions on $(L\cap H)\bb L$ is $L'$-invariant. From the properties of the constant term of smooth functions on $(L\cap H)\bb L$ and the characterization of the constant term of the elements of $V$ (\cite{D2} Proposition 3.14), one deduces the latter half of our first  claim. \\ As the  exponents of $(f_P)_{\vert M\cap L'}$ are the restrictions to $A_M\cap L'$ of exponents of $f_P$,  one deduces the second part of our claim. This achieves to prove (\ref{fPLL}). \\
Together with (\ref{ll'}) this implies
\ber\label{FLL}  If   $F$ be  a family of type I of tempered functions on $(H\cap L)\bb L$ then $F_{\vert L'}$ is a family of type I of tempered functions on $(H\cap L')\bb L'$.  \eer 
Let us prove that there exists constant $C, C'>0$  and $d\in \N$ such that 
\ber \label{nll} $$C N_d((H\cap L )l)^{-1} \Theta_L ((H\cap L)l) \leq \Theta_{L'} ((H\cap L' )l)$$
$$\leq C' N_d((H\cap L)l) \Theta_L ((H\cap L)l),  l\in L' .$$\eer 
From (\ref{Theta}) applied to $L$ and $L'$ this holds for $l\in A^-_{P'}$ for every minimal $\si$-parabolic  subgroup  $P'$  of $L'$.
 From this fact, from  the Cartan decomposition and the invariance of $\Theta_{L'}$ (resp., $\Theta_L$) by a compact open subgroup  of $L'$ (resp., $L$) the required inequality is a consequence of the following assertion applied to $L$ and $L'$.
 \ber  There exists $d\in \N$, and for all $g_1\in G$ there exists $c,c'>0$ such that:
$$cN_{-d} (Hg) \Theta_G(Hg) \leq \Theta_G(Hgg_1)\leq c'\Theta_G(g) N_d (Hg), g \in G.$$ \eer 
The right inequality is simply (\ref{lem3i}). The left inequality follows from the right one applied to $ gg_1$ instead of $g$ and $g_1^{-1}$ instead of $g_1$ as $N(Hgg_1)\leq CN(Hg) N(Hg_1)$ for some $C>0$. \\
% Your bilbigraphy           %<-------------------

\vskip 0,3cm
\begin{tabular}{lll}
Patrick Delorme& \hspace{1cm}   &Pascale Harinck   \\    
Aix Marseille Universit\'e&  \hspace{1cm}   &Ecole polytechnique\\
CNRS-IML, FRE 3529& \hspace{1cm} &CNRS-CMLS, UMR 7640\\
163 Avenue de Luminy&  \hspace{1cm}    &Route de Saclay\\ 
13288 Marseille Cedex 09& \hspace{1cm}     &91128 Palaiseau Cedex\\  France.& \hspace{1cm}  &France.\\            delorme@iml.univ-mrs.fr            &  \hspace{1cm} &harinck@math.polytechnique.fr\end{tabular}

\end{document}